\documentclass[twoside,12pt]{article}
\usepackage{amssymb}
\usepackage{amsthm}
\usepackage{amsmath}
\usepackage{a4wide}
\usepackage[all]{xy}
\usepackage{fancyhdr}

\newtheorem{theorem}{Theorem}[section]
\newtheorem{proposition}[theorem]{Proposition}
\newtheorem{corollary}[theorem]{Corollary}
\newtheorem{lemma}[theorem]{Lemma}
\theoremstyle{definition}
\newtheorem*{notation}{Notation}
\newtheorem*{Beweis}{Proof}
\newtheorem{definition}[theorem]{Definition}
\newtheorem{punto}[theorem]{}
\theoremstyle{remark}
\newtheorem{remark}[theorem]{Remark}
\newtheorem{ex}[theorem]{Example}
\newtheorem{c-ex}[theorem]{Counterexample}
\newtheorem{c-exs}[theorem]{Counterexamples}
\newtheorem{remarks}[theorem]{Remarks}
\CompileMatrices  \swapnumbers
\CompileMatrices
\input{tcilatex}

\begin{document}

\title{Fully Coprime Comodules and Fully Coprime\ Corings\thanks{%
MSC (2000): 16W30, 16N60, 16A53\newline
Keywords: Fully Coprime (Fully Cosemiprime) Corings, Prime (Semiprime)
Corings, Fully Coprime (Fully Cosemiprime) Comodules, Prime (Semiprime)
Comodules, Fully Coprime Spectrum, Fully Coprime Coradical}}
\author{\textbf{Jawad Y. Abuhlail}\thanks{%
Supported by King Fahd University of Petroleum $\&$ Minerals, Research
Project $\#$ INT/296} \\
%EndAName
Department of Mathematical Sciences, Box $\#\ $5046\\
King Fahd University of Petroleum $\&$ Minerals\\
31261 Dhahran - Saudi Arabia\\
abuhlail@kfupm.edu.sa}
\date{}
\maketitle

\begin{abstract}
Prime objects were defined as generalization of simple objects in the
categories of rings (modules). In this paper we introduce and investigate
what turns out to be a suitable generalization of simple corings (simple
comodules), namely \emph{fully coprime corings} (\emph{fully coprime
comodules}). Moreover, we consider several \emph{primeness} notions in the
category of comodules of a given coring and investigate their relations with
the fully coprimeness and the simplicity of these comodules. These notions
are applied then to study primeness and coprimeness properties of a given
coring, considered as an object in its category of right (left) comodules.
\end{abstract}

\section{Introduction}

\qquad \emph{Prime ideals} play a central role in the Theory of Rings. In
particular, \emph{localization of commutative rings} at prime ideals is an
essential tool in Commutative Algebra. One goal of this paper is to
introduce a suitable dual notion of \emph{coprimeness }for corings over
arbitrary (not necessarily commutative) ground rings as a first step towards
developing a theory of \emph{colocalization of corings}, which seems till
now to be far from reach.

The classical notion of a prime ring was generalized, in different ways, to
introduce prime objects in the category of modules of a given ring (see \cite%
[Section 13]{Wis96}). A main goal of this paper is to introduce \emph{coprime%
} objects, which generalize simple objects, in the category of corings
(comodules). As there are several \emph{primeness} properties of modules of
a given ring, we are led as well to several \emph{primeness} and \emph{%
coprimeness} properties of comodules of a coring. We investigate these
different properties and clarify the relations between them.

\emph{Coprime coalgebras} over base fields were introduced by R. Nekooei and
L. Torkzadeh in \cite{NT01} as a generalization of simple coalgebras: simple
coalgebras are coprime;\ and finite dimensional coprime coalgebras are
simple. These coalgebras, which we call here \emph{fully coprime,} were
defined using the so called \emph{wedge product} of subcoalgebras and can be
seen as dual to prime algebras: a coalgebra $C$ over a base field is coprime
if and only if its dual algebra $C^{\ast }$ is prime. Coprime coalgebras
were considered also by P. Jara et. al. in their study of representation
theory of coalgebras and path coalgebras \cite{JMR}.

For a coring $\mathcal{C}$ over a QF ring $A$ such that $_{A}\mathcal{C}$ ($%
\mathcal{C}_{A}$) is projective, we observe in Proposition \ref{wedge=} that
if $K,L\subseteq \mathcal{C}$ are any $A$-subbimodules that are right (left)
$\mathcal{C}$-coideals as well and satisfy suitable purity conditions, then
the \emph{wedge product} $K\wedge L,$ in the sense of \cite{Swe69}, is
nothing but their \emph{internal coproduct }$(K:_{\mathcal{C}^{r}}L)$ ($(K:_{%
\mathcal{C}^{l}}L)$) in the category of right (left) $\mathcal{C}$%
-comodules, in the sense of \cite{RRW05}. This observation suggests
extending the notion of fully coprime coalgebras over base fields to \emph{%
fully coprime corings} over arbitrary ground rings by replacing the wedge
product of subcoalgebras with the internal coproduct of subbicomodules. We
also extend that notion to \emph{fully coprime comodules} using the internal
coproduct of fully invariant subcomodules. Using the internal coproduct of a
bicoideal of a coring (a fully invariant subcomodule of a comodule) with
itself enables us to introduce \emph{fully cosemiprime corings} (\emph{fully
cosemiprime comodules}). Dual to prime radicals of rings (modules), we
introduce and investigate the \emph{fully coprime coradicals} of corings
(comodules).

This article is divided as follows: after this first introductory section,
we give in the second section some definitions and recall some needed
results from the Theory of Rings and Modules as well as from the Theory of
Corings and Comodules. As a coalgebra $C$ over a base field is fully coprime
if and only if its dual algebra $C^{\ast }\simeq \mathrm{End}^{C}(C)^{op}$
is prime, see \cite[Proposition 1.2]{NT01}, we devote the third section to
the study of primeness properties of the ring of $\mathcal{C}$-colinear
endomorphisms $\mathrm{E}_{M}^{\mathcal{C}}:=\mathrm{End}^{\mathcal{C}%
}(M)^{op}$ of a given right $\mathcal{C}$-comodule $M$ of a coring $\mathcal{%
C}.$ Given a coring $\mathcal{C},$ we say a non-zero right $\mathcal{C}$%
-comodule $M$ is $\mathrm{E}$\emph{-prime} (respectively $\mathrm{E}$\emph{%
-semiprime, completely } $\mathrm{E}$\emph{-prime, completely }$\mathrm{E}$%
\emph{-semiprime}), provided the ring $\mathrm{E}_{M}^{\mathcal{C}}:=\mathrm{%
End}^{\mathcal{C}}(M)^{op}$ is prime (respectively semiprime, domain,
reduced). In case $M$ is \emph{self-cogenerator}, Theorem \ref{iff-M}
provides sufficient and necessary conditions for $M$ to be $\mathrm{E}$%
-prime (respectively $\mathrm{E}$-semiprime, completely $\mathrm{E}$-prime,
completely $\mathrm{E}$-semiprime). Under suitable conditions, we clarify in
Theorem \ref{E-coprime-irr} the relation between $\mathrm{E}$-prime and
irreducible comodules. In the fourth section we present and study \emph{%
fully coprime} (\emph{fully cosemiprime}) \emph{comodules} using the \emph{%
internal coproduct} of fully invariant subcomodules. Let $\mathcal{C}$ be a
coring and $M$ be a non-zero right $\mathcal{C}$-comodule. A fully invariant
non-zero $\mathcal{C}$-subcomodule $K\subseteq M$ will be called fully$%
\mathcal{\ }M$-coprime (fully $M$-cosemiprime), iff for any (equal) fully
invariant $\mathcal{C}$-subcomodules $X,Y\subseteq M$ with $K\subseteq
(X:_{M}^{\mathcal{C}}Y),$ we have $K\subseteq X$ or $K\subseteq Y,$ where $%
(X:_{M}^{\mathcal{C}}Y)$ is the internal coproduct of $X,Y$ in the category
of right $\mathcal{C}$-comodules. We call the non-zero right $\mathcal{C}$%
-comodule $M$ \emph{fully coprime }(\emph{fully cosemiprime}), iff $M$ is
fully$\mathcal{\ }M$-coprime (fully $M$-cosemiprime). The notion of fully
coprimeness (fully cosemiprimeness) in the category of left $\mathcal{C}$%
-comodules is defined analogously. Theorem \ref{coprime-tau} clarifies the
relation between fully coprime (fully cosemiprime) and $\mathrm{E}$-prime ($%
\mathrm{E}$-semiprime) comodules under suitable conditions. We define the
\emph{fully coprime spectrum} of $M$ as the class of all fully$\mathcal{\ }M$%
-coprime $\mathcal{C}$-subcomodules of $M$ and the \emph{fully coprime
coradical }of $M$ as the sum of all fully$\mathcal{\ }M$-coprime $\mathcal{C}
$-subcomodules. In Proposition \ref{Prad=CPcorad} we clarify the relation
between the fully coprime coradical of $M$ and the prime radical of $\mathrm{%
E}_{M}^{\mathcal{C}},$ in case $M$ is intrinsically injective
self-cogenerator and $\mathrm{E}_{M}^{\mathcal{C}}$ is right Noetherian.
Fully coprime comodules turn to be a generalization of simple comodules:
simple comodules are trivially fully coprime; and Theorem \ref%
{coprime=simple} (2) shows that if the ground ring $A$ is right Artinian and
$_{A}\mathcal{C}$ is locally projective, then a non-zero finitely generated
self-injective self-cogenerator right $\mathcal{C}$-comodule $M$ is fully
coprime if and only if $M$ is simple as a $(^{\ast }\mathcal{C},\mathrm{E}%
_{M}^{\mathcal{C}})$-bimodule. Under suitable conditions, we clarify in
Theorem \ref{coprime-irr} the relation between fully coprime and irreducible
comodules. In the fifth section we introduce and study several primeness and
coprimeness properties of a non-zero coring $\mathcal{C},$ considered as an
object in the category $\mathbb{M}^{\mathcal{C}}$ of right $\mathcal{C}$%
-comodules and as an object in the category $^{\mathcal{C}}\mathbb{M}$ of
left $\mathcal{C}$-comodules. We define the internal coproducts of $\mathcal{%
C}$-bicoideals, i.e. $(\mathcal{C},\mathcal{C})$-subbicomodules of $\mathcal{%
C},$ in $\mathbb{M}^{\mathcal{C}}$ and in $^{\mathcal{C}}\mathbb{M}$ and use
them to introduce the notions of \emph{fully coprime} (\emph{fully
cosemiprime}) $\mathcal{C}$\emph{-bicoideals} and \emph{fully coprime} (%
\emph{fully cosemiprime}) \emph{corings}. Moreover, we introduce and study
the \emph{fully coprime spectrum }and the \emph{fully coprime coradical }of $%
\mathcal{C}$ in $\mathbb{M}^{\mathcal{C}}$ (in $^{\mathcal{C}}\mathbb{M})$
and clarify their relations with the prime spectrum and the prime radical of
$\mathcal{C}^{\ast }$ ($^{\ast }\mathcal{C}$). We investigate several
coprimeness (cosemiprimeness) and primeness (semiprimeness) notions for $%
\mathcal{C}$ and clarify their relations with the simplicity
(semisimplicity) of the coring under consideration. In Theorems \ref{sn-cor}
we give sufficient and necessary conditions for the dual ring $\mathcal{C}%
^{\ast }$ ($^{\ast }\mathcal{C}$) of $\mathcal{C\ }$to be prime
(respectively semiprime, domain, reduced). In case the ground ring $A$ is a
QF ring, $_{A}\mathcal{C},\mathcal{C}_{A}$ are locally projective and $%
\mathcal{C}^{\ast }$ is right Artinian, $^{\ast }\mathcal{C}$ is left
Artinian, we show in Theorem \ref{A-C*-Art} that $\mathcal{C}^{r}$ is fully
coprime if and only if $\mathcal{C}$ is simple if and only if $\mathcal{C}%
^{l}$ is fully coprime.

Throughout, $R$ is a commutative ring with $1_{R}\neq 0_{R},$ $A$ is an
arbitrary but fixed unital $R$-algebra and $\mathcal{C}$ is a non-zero $A$%
-coring. All rings have unities preserved by morphisms of rings and all
modules are unital. Let $T$ be a ring and denote with $_{T}\mathbb{M}$ ($%
\mathbb{M}_{T}$) the category of left (right) $T$-modules. For a left
(right) $T$-module $M,$ we denote with $\sigma \lbrack _{T}M]\subseteq $ $%
_{T}\mathbb{M}$ ($\sigma \lbrack M_{T}\subseteq \mathbb{M}_{T}]$) the \emph{%
full }subcategory of $M$-\emph{subgenerated} left (right) $T$-modules; see
\cite{Wis91} and \cite{Wis96}.

\section{Preliminaries}

In this section we introduce some definitions, remarks and lemmas to which
we refer later.

\begin{punto}
(\cite{Z-H76}) An $A$-module $W$ is called \emph{locally projective} (in the
sense of B. Zimmermann-Huisgen \cite{Z-H76}), if for every diagram
\begin{equation*}
\xymatrix{0 \ar[r] & F \ar@{.>}[dr]_{g' \circ \iota} \ar[r]^{\iota} & W
\ar[dr]^{g} \ar@{.>}[d]^{g'} & & \\ & & L \ar[r]_{\pi} & N \ar[r] & 0}
\end{equation*}%
with exact rows and $F$ f.g.: for every $A$-linear map $g:W\rightarrow N,$
there exists an $A$-linear map $g^{\prime }:W\rightarrow L$, such that the
entstanding parallelogram is commutative. Note that every projective $A$%
-module is locally projective. Moreover, every locally projective $A$-module
is flat and $A$-cogenerated.
\end{punto}

\subsection*{Prime and coprime modules}

\begin{definition}
Let $T$ be a ring. A proper ideal $P\vartriangleleft T$ is called

\emph{prime}, iff for any two ideals $I,J\vartriangleleft T$ with $%
IJ\subseteq P,$ either $I\subseteq P$ or $J\subseteq P;$

\emph{semiprime, }iff for any ideal $I\vartriangleleft T$ with $%
I^{2}\subseteq P,$ we have $I\subseteq P;$

\emph{completely prime}, iff for any $f,g\in P$ with $fg\in P,$ either $f\in
P$ or $g\in P;$

\emph{completely semiprime}, iff for any $f\in T$ with $f^{2}\in P,$ we have
$f\in P.$

The ring $T$ is called \emph{prime} (respectively \emph{semiprime}, \emph{%
domain}, \emph{reduced}), iff $(0_{T})\vartriangleleft T$ is prime
(respectively semiprime, completely prime, completely semiprime).
\end{definition}

\begin{punto}
Let $T$ be a ring. With $\mathrm{Max}(T)$ (resp. $\mathrm{Max}_{r}(T),$ $%
\mathrm{Max}_{l}(T)$) we denote the class of maximal two-sided $T$-ideals
(resp. maximal right $T$-ideals, maximal left $T$-ideals) and with $\mathrm{%
Sepc}(T)$ the prime spectrum of $T$ consisting of all prime ideals of $T.$
The \emph{Jacobson radical }of $T$ is denoted by $\mathrm{Jac}(T)$ and the
\emph{prime radical} of $T$ by $\mathrm{Prad}(T).$ Notice that the ring $T$
is semiprime if and only if $\mathrm{Prad}(T)=0.$
\end{punto}

There are various notions of prime and coprime modules in the literature;
see \cite[Section 13]{Wis96} for more details. In this paper we adopt the
notion of \emph{prime modules} due to R. Johnson \cite{Joh53} and its dual
notion of \emph{coprime modules} considered recently by S. Annin \cite{Ann}.

\begin{definition}
Let $T$ be a ring. A non-zero $T$-module $M$ will be called

\emph{prime,} iff $\mathrm{ann}_{T}(K)=\mathrm{ann}_{T}(M)$ for every
non-zero $T$-submodule $0\neq K\subseteq M;$

\emph{coprime,} iff $\mathrm{ann}_{T}(M/K)=\mathrm{ann}_{T}(M)$ for every
proper $T$-submodule $K\varsubsetneqq M;$

\emph{diprime,} iff $\mathrm{ann}_{T}(K)=\mathrm{ann}_{T}(M)$ or $\mathrm{ann%
}_{T}(M/K)=\mathrm{ann}_{T}(M)$ for every non-trivial $T$-submodule $0\neq
K\varsubsetneqq M;$

\emph{strongly prime,} iff $M\in \sigma \lbrack K]$ for every non-zero $T$%
-submodule $0\neq K\subseteq M;$

\emph{semiprime}, iff $M/\mathcal{T}_{K}(M)\in \sigma \lbrack K]$ for every
cyclic $T$-submodule $K\subseteq M;$

\emph{strongly semiprime},\emph{\ }iff $M/\mathcal{T}_{K}(M)\in \sigma
\lbrack K]$ for every $T$-submodule $K\subseteq M.$
\end{definition}

It's well known that for every prime (coprime) $T$-module $M,$ the
associated quotient ring $\overline{T}:=T/\mathrm{ann}_{T}(M)$ is prime. In
fact we have more:

\begin{proposition}
\label{Tbar}\emph{(\cite[Proposition 1.1]{Lom05})} Let $T$ be a ring and $M$
be a non-zero $T$-module. Then the following are equivalent:

\begin{enumerate}
\item $\overline{T}:=T/\mathrm{ann}_{T}(M)$ is a prime ring;

\item $M$ is diprime;

\item For every fully invariant $T$-submodule $K\subseteq M$ that is $M$%
-generated as an $\mathrm{End}_{T}(M)$-module, $\mathrm{ann}_{T}(K)=\mathrm{%
ann}_{T}(M)$ or $\mathrm{ann}_{T}(M/K)=\mathrm{ann}_{T}(M).$
\end{enumerate}
\end{proposition}

\begin{remark}
\label{*-**}Let $T$ be a ring and consider the following conditions for a
non-zero $T$-module $M:$%
\begin{equation*}
\mathrm{ann}_{T}(M/K)\neq \mathrm{ann}_{T}(M)\text{ for every non-zero }T%
\text{-submodule }0\neq K\subseteq M\text{ (*)}
\end{equation*}%
\begin{equation*}
\mathrm{ann}_{T}(K)\neq \mathrm{ann}_{T}(M)\text{ for every proper }T\text{%
-submodule }K\varsubsetneqq M\text{ (**).}
\end{equation*}%
We introduce condition (**) as dual to condition (*), which is due to
Wisbauer \cite[Section 13]{Wis96}. Modules satisfying either of these
conditions allow further conclusions from the primeness (coprimeness)
properties: by Proposition \ref{Tbar}, a $T$-module $M$ satisfying condition
(*)\emph{\ }(condition (**)) is prime (coprime) if and only if $\overline{T}%
:=T/\mathrm{ann}_{T}(M)$ is prime.
\end{remark}

\subsection*{Corings and comodules}

Fix a non-zero $A$-coring $(\mathcal{C},\Delta ,\varepsilon ).$ With $%
\mathbb{M}^{\mathcal{C}}$ ($^{\mathcal{C}}\mathbb{M}$) we denote the
category of right (left) $\mathcal{C}$-comodules with the $\mathcal{C}$%
-colinear morphisms and by $\mathcal{C}^{r}$ ($\mathcal{C}^{l}$) we mean the
coring $\mathcal{C},$ considered as an object in $\mathbb{M}^{\mathcal{C}}$ (%
$^{\mathcal{C}}\mathbb{M}$). For a right (left) $\mathcal{C}$-comodule $M$
we denote with $\mathrm{E}_{M}^{\mathcal{C}}:=\mathrm{End}^{\mathcal{C}%
}(M)^{op}$ ($_{M}^{\mathcal{C}}\mathrm{E}:=$ $^{\mathcal{C}}\mathrm{End}(M)$%
) the ring of all $\mathcal{C}$-colinear endomorphisms of $M$ with
multiplication the opposite (usual) composition of maps and call an $R$%
-submodule $X\subseteq M$ \emph{fully invariant,} iff $f(X)\subseteq X$ for
every $f\in \mathrm{E}_{M}^{\mathcal{C}}$ ($f\in $ $_{M}^{\mathcal{C}}%
\mathrm{E}$).

In module categories, monomorphisms are injective maps. In comodule
categories this is not the case in general. In fact we have:

\begin{remark}
\label{mono}For any coring $\mathcal{C}$ over a ground ring $A,$ the module $%
_{A}\mathcal{C}$ is flat if and only if every monomorphism in $\mathbb{M}^{%
\mathcal{C}}$ is injective (e.g. \cite[Proposition 1.10]{Abu03}). In this
case, $\mathbb{M}^{\mathcal{C}}$ is a Grothendieck category with kernels
formed in the category of right $A$-modules and given a right $\mathcal{C}$%
-comodule $M,$ the intersection $\bigcap_{\lambda \in \Lambda }M_{\lambda
}\subseteq M$ of any family $\{M_{\lambda }\}_{\Lambda }$ of $\mathcal{C}$%
-subcomodules of $M$ is again a $\mathcal{C}$-subcomodule.
\end{remark}

\begin{definition}
Let $_{A}\mathcal{C}$ ($\mathcal{C}_{A}$) be flat. We call a non-zero right
(left) $\mathcal{C}$-subcomodule $M$

\emph{simple, }iff $M$ has no non-trivial $\mathcal{C}$-subcomodules;

\emph{semisimple}, iff $M=\mathrm{Soc}(M)$ where%
\begin{equation}
\mathrm{Soc}(M):=\bigoplus \{K\subseteq M\mid K\text{ is a simple }\mathcal{C%
}\text{-subcomodule}\}.  \label{soc(M)}
\end{equation}

The right (left) $\mathcal{C}$-subcomodule $\mathrm{Soc}(M)\subseteq M$
defined in (\ref{soc(M)}) is called the \emph{socle} of $M.$We call a
non-zero right (left) $\mathcal{C}$-subcomodule $0\neq K\subseteq M$ \emph{%
essential in }$M,$ and write $K\vartriangleleft _{e}M,$ provided $K\cap
\mathrm{Soc}(M)\neq 0.$
\end{definition}

\begin{lemma}
\label{MotC-inj}\emph{(\cite[Proposition 1.10]{Abu03})} If $A_{A}$ is
injective \emph{(}cogenerator\emph{)} and $N$ is a right $A$-module, then
the canonical right $\mathcal{C}$-comodule $M:=(N\otimes _{A}\mathcal{C}%
,id\otimes _{A}\Delta _{\mathcal{C}})$ is injective \emph{(}cogenerator\emph{%
)} in $\mathbb{M}^{\mathcal{C}}.$ In particular, if $A_{A}$ is injective
\emph{(}cogenerator\emph{)} then $\mathcal{C}\simeq A\otimes _{A}\mathcal{C}$
is injective \emph{(}cogenerator\emph{)} in $\mathbb{M}^{\mathcal{C}}.$
\end{lemma}

For an $A$-coring $\mathcal{C},$ the dual module $^{\ast }\mathcal{C}:=%
\mathrm{Hom}_{A-}(\mathcal{C},A)$ ($\mathcal{C}^{\ast }:=\mathrm{Hom}_{-A}(%
\mathcal{C},A)$) of left (right) $A$-linear maps from $\mathcal{C}$ to $A$
is a ring under the so called \emph{convolution product}. We remark here
that the multiplications used below are opposite to those in previous papers
of the author, e.g. \cite{Abu03}, and are consistent with the ones in \cite%
{BW03}.

\begin{punto}
\textbf{Dual rings of corings. }Let $(\mathcal{C},\Delta ,\varepsilon )$ be
an $A$-coring. Then $^{\ast }\mathcal{C}:=\mathrm{Hom}_{A-}(\mathcal{C},A)$
(respectively $\mathcal{C}^{\ast }:=\mathrm{Hom}_{-A}(\mathcal{C},A)$) is an
$A^{op}$-ring with multiplication
\begin{equation*}
(f\ast ^{l}g)(c)=\sum f(c_{1}g(c_{2}))\text{ (respectively }(f\ast
^{r}g)(c)=\sum g(f(c_{1})c_{2})
\end{equation*}%
and unity $\varepsilon .$ The coring $\mathcal{C}$ is a $(^{\ast }\mathcal{C}%
,\mathcal{C}^{\ast })$-bimodule through the left $^{\ast }\mathcal{C}$%
-action (respectively the right $\mathcal{C}^{\ast }$-action):%
\begin{equation*}
f\rightharpoonup c:=\sum c_{1}f(c_{2})\text{ for all }f\in \text{ }^{\ast }%
\mathcal{C}\text{ (respectively }c\leftharpoonup g:=\sum g(c_{1})c_{2}\text{
for all }g\in \mathcal{C}^{\ast }\text{).}
\end{equation*}
\end{punto}

\begin{punto}
Let $M$ be a right (left) $\mathcal{C}$-comodule. Then $M$ is a left $^{\ast
}\mathcal{C}$-module (a right $\mathcal{C}^{\ast }$-module) under the left
(right) action
\begin{equation*}
f\rightharpoonup m:=\sum m_{<0>}f(m_{<1>})\text{ for all }f\in \text{ }%
^{\ast }\mathcal{C}\text{ \ (}m\leftharpoonup g:=\sum g(m_{<-1>})m_{<0>}%
\text{ for all }g\in \mathcal{C}^{\ast }\text{).}
\end{equation*}%
Notice that $M$ is a $(^{\ast }\mathcal{C},\mathrm{E}_{M}^{\mathcal{C}})$%
-bimodule (a $(\mathcal{C}^{\ast },$ $_{M}^{\mathcal{C}}\mathrm{E})$%
-bimodule) in the canonical way. A right (left) $\mathcal{C}$-subcomodule $%
K\subseteq M$ is said to be \emph{fully invariant}, provided $K$ is a $%
(^{\ast }\mathcal{C},\mathrm{E}_{M}^{\mathcal{C}})$-subbimodule ($(\mathcal{C%
}^{\ast },$ $_{M}^{\mathcal{C}}\mathrm{E})$-subbimodule) of $M.$ Since $%
\mathbb{M}^{\mathcal{C}}$ ($^{\mathcal{C}}\mathbb{M}$) has cokernels, we
conclude that for any $f\in \mathrm{E}_{M}^{\mathcal{C}}$ (any $g\in $ $%
_{M}^{\mathcal{C}}\mathrm{E}$), $Mf:=f(M)\subseteq M$ ($gM:=g(M)\subseteq M$%
) is a right (left) $\mathcal{C}$-subcomodule and that for any right ideal $%
I\vartriangleleft _{r}\mathrm{E}_{M}^{\mathcal{C}}$ (left ideal $%
J\vartriangleleft _{l}$ $_{M}^{\mathcal{C}}\mathrm{E}$) we have a
fully-invariant right (left) $\mathcal{C}$-subcomodule $MI\subseteq M$ ($%
JM\subseteq M$).
\end{punto}

\begin{proposition}
\label{comod=sg}\emph{(\cite[Theorems 2.9, 2.11]{Abu03}) }For any $A$-coring
$\mathcal{C}$ we have

\begin{enumerate}
\item $\mathbb{M}^{\mathcal{C}}\simeq \sigma \lbrack \mathcal{C}_{^{\ast }%
\mathcal{C}^{op}}]\simeq \sigma \lbrack _{^{\ast }\mathcal{C}}\mathcal{C}]$
if and only if $_{A}\mathcal{C}$ is locally projective.

\item $^{\mathcal{C}}\mathbb{M}\simeq \sigma \lbrack _{\mathcal{C}^{\ast op}}%
\mathcal{C}]\simeq \sigma \lbrack \mathcal{C}_{\mathcal{C}^{\ast }}]$ if and
only if $\mathcal{C}_{A}$ is locally projective.
\end{enumerate}
\end{proposition}

\begin{notation}
Let $M$ be a right $\mathcal{C}$-comodule. We denote with $\mathcal{C}(M)$ ($%
\mathcal{C}_{f.i.}(M)$) the class of (fully invariant) $\mathcal{C}$%
-subcomodules of $M$ and with $\mathcal{I}_{r}(\mathrm{E}_{M}^{\mathcal{C}})$
($\mathcal{I}_{t.s.}(\mathrm{E}_{M}^{\mathcal{C}})$) the class of right
(two-sided) ideals of $\mathrm{E}_{M}^{\mathcal{C}}.$ For $\varnothing \neq
K\subseteq M,$ $\varnothing \neq I\subseteq \mathrm{E}_{M}^{\mathcal{C}}$
set
\begin{equation*}
\mathrm{An}(K):=\{f\in \mathrm{E}_{M}^{\mathcal{C}}\mid f(K)=0\},\text{ }%
\mathrm{Ke}(I):=\bigcap \{\mathrm{Ker}(f)\mid f\in I\}.
\end{equation*}
\end{notation}

The following notions for right $\mathcal{C}$-comodules will be used in the
sequel. The analogous notions for left $\mathcal{C}$-comodules can be
defined analogously:

\begin{definition}
\label{injectivity}Let $_{A}\mathcal{C}$ be flat. We say a right $\mathcal{C}
$-comodule $M$ is

\emph{self-injective, }iff for every $\mathcal{C}$-subcomodule $K\subseteq
M, $ every $\mathcal{C}$-colinear morphism $f\in \mathrm{Hom}^{\mathcal{C}%
}(K,M) $ extends to a $\mathcal{C}$-colinear endomorphism $\widetilde{f}\in
\mathrm{End}\mathcal{^{C}}(M);$

\emph{semi-injective, }iff for every monomorphism $0\longrightarrow N\overset%
{h}{\longrightarrow }M$ in $\mathbb{M}^{\mathcal{C}},$ where $N$ is a factor
$\mathcal{C}$-comodule of $M,$ and every $f\in \mathrm{Hom}^{\mathcal{C}%
}(N,M),$ $\exists $ $\widetilde{f}\in \mathrm{End}\mathcal{^{C}}(M)$ such
that $\widetilde{f}\circ h=f;$

\emph{self-projective}, iff for every $\mathcal{C}$-subcomodule $K\subseteq
M,$ and every $g\in \mathrm{Hom}^{\mathcal{C}}(M,M/K),$ $\exists $ $%
\widetilde{g}\in \mathrm{End}\mathcal{^{C}}(M)$ such that $\pi _{K}\circ
\widetilde{g}=g;$

\emph{self-cogenerator,} iff $M$ cogenerates all of its factor $\mathcal{C}$%
-comodules;

\emph{self-generator,} iff $M$ generates each of its $\mathcal{C}$%
-subcomodules;

\emph{coretractable,} iff $\mathrm{Hom}^{\mathcal{C}}(M/K,M)\neq 0$ for
every proper $\mathcal{C}$-subcomodule $K\subsetneqq M;$

\emph{retractable, }iff $\mathrm{Hom}^{\mathcal{C}}(M,K)\neq 0$ for every
non-zero $\mathcal{C}$-subcomodule $0\neq K\subseteq M;$

\emph{intrinsically injective}, iff $\mathrm{AnKe}(I)=I$ for every f.g.
right ideal $I\vartriangleleft \mathrm{E}_{M}^{\mathcal{C}}.$
\end{definition}

\qquad The following result follows immediately from (\cite[31.11, 31.12]%
{Wis91}) and Proposition \ref{comod=sg}:

\begin{proposition}
\label{End} Let $_{A}\mathcal{C}$ be locally projective, $M$ be a non-zero
right $\mathcal{C}$-comodule and consider the ring $\mathrm{E}_{M}^{\mathcal{%
C}}:=\mathrm{End}^{\mathcal{C}}(M)^{op}=\mathrm{End}(_{^{\ast }\mathcal{C}%
}M)^{op}.$

\begin{enumerate}
\item If $M$ is Artinian and self-injective, then $\mathrm{E}_{M}^{\mathcal{C%
}}$ is right Noetherian.

\item If $M$ is Artinian, self-injective and self-projective, then $\mathrm{E%
}_{M}^{\mathcal{C}}$ is right Artinian.

\item If $M$ is semi-injective and satisfies the ascending chain condition
for annihilator $\mathcal{C}$-subcomodules, then $\mathrm{E}_{M}^{\mathcal{C}%
}$ is semiprimary.
\end{enumerate}
\end{proposition}

\subsection*{Annihilator conditions for comodules\qquad}

\qquad Analogous to the \emph{annihilator conditions} for modules (e.g. \cite%
[28.1]{Wis91}), the following result gives some annihilator conditions for
comodules.

\begin{punto}
\label{An-Ke} Let $_{A}\mathcal{C}$ be flat, $M$ be a right $\mathcal{C}$%
-comodule and consider the order-reversing mappings
\begin{equation}
\mathrm{An}(-):\mathcal{C}(M)\rightarrow \mathcal{I}_{r}(\mathrm{E}_{M}^{%
\mathcal{C}})\text{ and }\mathrm{Ke}(-):\mathcal{I}_{r}(\mathrm{E}_{M}^{%
\mathcal{C}})\rightarrow \mathcal{C}(M).  \label{An-map-Ke}
\end{equation}

\begin{enumerate}
\item For every $K\in \mathcal{C}_{f.i.}(M)$ ($I\in \mathcal{I}_{t.s.}(%
\mathrm{E}_{M}^{\mathcal{C}})$), we have $\mathrm{An}(K)\in \mathcal{I}%
_{t.s.}(\mathrm{E}_{M}^{\mathcal{C}})$ ($\mathrm{Ke}(I)\in \mathcal{C}%
_{f.i.}(M)$). Moreover $\mathrm{An}(-)$ and $\mathrm{Ke}(-)$ induce
bijections
\begin{equation*}
\begin{tabular}{lllll}
$\mathcal{A}(\mathrm{E}_{M}^{\mathcal{C}})$ & $:=\{\mathrm{An}(K)|\text{ }%
K\in \mathcal{C}(M)\}$ & $\leftrightarrow $ & $\mathcal{K}(M)$ & $:=\{%
\mathrm{Ke}(I)|$ $I\in \mathcal{I}_{r}(\mathrm{E}_{M}^{\mathcal{C}})\};$ \\
$\mathcal{A}_{t.s.}(\mathrm{E}_{M}^{\mathcal{C}})$ & $:=\{\mathrm{An}(K)|%
\text{ }K\in \mathcal{C}_{f.i.}(M)\}$ & $\leftrightarrow $ & $\mathcal{K}%
_{f.i.}(M)$ & $:=\{\mathrm{Ke}(I)|\text{ }I\in \mathcal{I}_{t.s.}(\mathrm{E}%
_{M}^{\mathcal{C}})\}.$%
\end{tabular}%
\end{equation*}

\item For any $\mathcal{C}$-subcomodule $K\subseteq M$ we have
\begin{equation*}
\mathrm{KeAn}(K)=K\text{ if and only if }M/K\text{ is }M\text{-cogenerated.}
\end{equation*}

\item If $M$ is self-injective, then

\begin{enumerate}
\item $\mathrm{An}(\bigcap\limits_{i=1}^{n}K_{i})=\sum\limits_{i=1}^{n}%
\mathrm{An}(K_{i})$ for any finite set of $\mathcal{C}$-subcomodules $%
K_{1},...,K_{n}\subseteq M.$

\item $M$ is intrinsically injective.
\end{enumerate}
\end{enumerate}
\end{punto}

\begin{remarks}
\label{inj-cor}let $_{A}\mathcal{C}$ be flat and $M$ be a right $\mathcal{C}$%
-comodule.

\begin{enumerate}
\item If $M$ is self-injective (self-cogenerator), then every \emph{fully
invariant} $\mathcal{C}$-subcomodule of $M$ is also self-injective
(self-cogenerator).

\item If $M$ is self-injective, then $M$ is semi-injective. If $M$ is
self-generator (self-cogenerator), then it is obviously retractable
(coretractable).

\item If $M$ is self-cogenerator ($M$ is intrinsically injective and $%
\mathrm{E}_{M}^{\mathcal{C}}$ is right Noetherian), then the mapping
\begin{equation*}
\mathrm{An}(-):\mathcal{C}(M)\rightarrow \mathcal{I}_{r}(\mathrm{E}_{M}^{%
\mathcal{C}})\text{ (}\mathrm{Ke}(-):\mathcal{I}_{r}(\mathrm{E}_{M}^{%
\mathcal{C}})\rightarrow \mathcal{C}(M)\text{)}
\end{equation*}%
is injective.

\item Let $M$ be self-injective. If $H\subsetneqq K\subseteq M$ are $%
\mathcal{C}$-subcomodules with $K$ coretractable and fully invariant in $M,$
then $\mathrm{An}(K)\subsetneqq \mathrm{An}(H):$ since $M$ is self-injective
and $K\subseteq M$ is fully invariant, we have a surjective morphism of $R$%
-algebras $\mathrm{E}_{M}^{\mathcal{C}}\rightarrow \mathrm{E}_{K}^{\mathcal{C%
}}\rightarrow 0,$ $f\mapsto f_{|_{K}},$ which induces a bijection $\mathrm{An%
}(H)/\mathrm{An}(K)\longleftrightarrow \mathrm{An}_{\mathrm{E}_{K}^{\mathcal{%
C}}}(H)\simeq \mathrm{Hom}^{\mathcal{C}}(K/H,K)\neq 0.$
\end{enumerate}
\end{remarks}

\section{$\mathrm{E}$-prime ($\mathrm{E}$-semiprime)\ Comodules}

\qquad In this section we study and characterize non-zero comodules, for
which the ring of colinear endomorphisms is prime (respectively semiprime,
domain, reduced). Throughout, we assume $\mathcal{C}$ is a non-zero $A$%
-coring with $_{A}\mathcal{C}$ flat, $M$ is a non-zero right $\mathcal{C}$%
-comodule and $\mathrm{E}_{M}^{\mathcal{C}}:=\mathrm{End}^{\mathcal{C}%
}(M)^{op}$ is the ring of $\mathcal{C}$-colinear endomorphisms of $M$ with
the opposite composition of maps. We remark that analogous results to those
obtained in this section can be obtained for left $\mathcal{C}$-comodules,
by symmetry.

\begin{definition}
We define a fully invariant non-zero $\mathcal{C}$-subcomodule $0\neq
K\subseteq M$ to be

$\mathrm{E}$\emph{-prime in }$M,$ iff $\mathrm{An}(K)\vartriangleleft
\mathrm{E}_{M}^{\mathcal{C}}$ is prime;

$\mathrm{E}$\emph{-semiprime in }$M,$ iff $\mathrm{An}(K)\vartriangleleft
\mathrm{E}_{M}^{\mathcal{C}}$ is semiprime;

\emph{completely }$\mathrm{E}$\emph{-prime in }$M,$ iff $\mathrm{An}%
(K)\vartriangleleft \mathrm{E}_{M}^{\mathcal{C}}$ is completely prime;

\emph{completely }$\mathrm{E}$\emph{-semiprime in }$M$, iff $\mathrm{An}%
(K)\vartriangleleft \mathrm{E}_{M}^{\mathcal{C}}$ is completely semiprime.
\end{definition}

\begin{definition}
We call the right $\mathcal{C}$-comodule $M$ $\mathrm{E}$\emph{-prime}
(respectively $\mathrm{E}$\emph{-semiprime}, \emph{completely }$\mathrm{E}$%
\emph{-prime}, \emph{completely }$\mathrm{E}$\emph{-semiprime}), provided $M$
is $\mathrm{E}$-prime in $M$ (respectively $\mathrm{E}$-semiprime in $M$,
completely $\mathrm{E}$-prime in $M$, completely $\mathrm{E}$-semiprime in $%
M $), equivalently iff $R$-algebra $\mathrm{E}_{M}^{\mathcal{C}}$ is prime
(respectively semiprime, domain, reduced).
\end{definition}

\begin{notation}
For the right $\mathcal{C}$-comodule $M$ we denote with \textrm{EP}$(M)$
(resp. \textrm{ESP}$(M),$ \textrm{CEP}$(M),$ \textrm{CESP}$(M)$) the class
of fully invariant $\mathcal{C}$-subcomodules of $M$ whose annihilator in $%
\mathrm{E}_{M}^{\mathcal{C}}$ is prime (resp. semiprime, completely prime,
completely semiprime).
\end{notation}

\begin{ex}
\label{ex-AnKe}Let $P\vartriangleleft \mathrm{E}_{M}^{\mathcal{C}}$ be a
proper two-sided ideal with $P=\mathrm{AnKe}(P)$ (e.g. $M$ intrinsically
injective and $P_{\mathrm{E}_{M}^{\mathcal{C}}}$ finitely generated) and
consider the fully invariant $\mathcal{C}$-subcomodule $0\neq K:=\mathrm{Ke}%
(P)\subseteq M.$ Assume $P\vartriangleleft \mathrm{E}_{M}^{\mathcal{C}}$ to
be prime (respectively semiprime, completely prime, completely semiprime).
Then $K\in \mathrm{EP}(M)$ (respectively $K\in \mathrm{ESP}(M),$ $K\in
\mathrm{CEP}(M),$ $K\in \mathrm{CESP}(M)$). If moreover $M$ is
self-injective, then we have isomorphisms of $R$-algebras%
\begin{equation*}
\mathrm{E}_{K}^{\mathcal{C}}\simeq \mathrm{E}_{M}^{\mathcal{C}}/\mathrm{An}%
(K)=\mathrm{E}_{M}^{\mathcal{C}}/\mathrm{AnKe}(P)=\mathrm{E}_{M}^{\mathcal{C}%
}/P,
\end{equation*}%
hence $K$ is $\mathrm{E}$-prime (respectively $\mathrm{E}$-semiprime,
completely $\mathrm{E}$-prime, completely $\mathrm{E}$-semiprime).
\end{ex}

For the right $\mathcal{C}$-comodule $M$ we have
\begin{equation}
\mathrm{CEP}(M)\subseteq \mathrm{EP}(M)\subseteq \mathrm{ESP}(M)\text{ and }%
\mathrm{CEP}(M)\subseteq \mathrm{CESP}(M)\subseteq \mathrm{ESP}(M).
\label{incl}
\end{equation}

\begin{remark}
The idea of Example \ref{ex-AnKe} can be used to construct \emph{%
counterexamples}, which show that the inclusions in (\ref{incl}) are \emph{%
in general} strict.
\end{remark}

\subsection*{The $\mathrm{E}$-prime coradical}

\begin{definition}
We define the $\mathrm{E}$\emph{-prime coradical} of the right $\mathcal{C}$%
-comodule $M$ as%
\begin{equation*}
\mathrm{EPcorad}(M)=\sum_{K\in \mathrm{EP}(M)}K.
\end{equation*}
\end{definition}

\begin{proposition}
\label{P-rad}Let $M$ be intrinsically injective. If $\mathrm{E}_{M}^{%
\mathcal{C}}$ is right Noetherian, then%
\begin{equation}
\mathrm{Prad}(\mathrm{E}_{M}^{\mathcal{C}})=\mathrm{An}(\mathrm{EPcorad}(M)).
\label{Prad=AnEP}
\end{equation}%
If moreover $M$ is self-cogenerator, then
\begin{equation}
\mathrm{EPcorad}(M)=\mathrm{Ke}(\mathrm{Prad}(\mathrm{E}_{M}^{\mathcal{C}})).
\label{EP=KePrad}
\end{equation}
\end{proposition}

\begin{Beweis}
If $K\in \mathrm{EP}(M),$ then $\mathrm{An}(K)\vartriangleleft \mathrm{E}%
_{M}^{\mathcal{C}}$ is a prime ideal (by definition). On the otherhand, if $%
P\vartriangleleft \mathrm{E}_{M}^{\mathcal{C}}$ is a prime ideal then $P=%
\mathrm{AnKe}(P)$ (since $P_{\mathrm{E}_{M}^{\mathcal{C}}}$ is finitely
generated and $M$ is intrinsically injective) and so $K:=\mathrm{Ke}(P)\in
\mathrm{EP}(M).$ It follows then that
\begin{equation*}
\begin{tabular}{lllll}
$\mathrm{Prad}(\mathrm{E}_{M}^{\mathcal{C}})$ & $=$ & $\bigcap_{P\in \mathrm{%
Spec}(\mathrm{E}_{M}^{\mathcal{C}})}P$ & $=$ & $\bigcap_{P\in \mathrm{Spec}(%
\mathrm{E}_{M}^{\mathcal{C}})}\mathrm{AnKe}(P)$ \\
& $=$ & $\bigcap_{K\in \mathrm{EP}(M)}\mathrm{AnKeAn}(K)$ & $=$ & $%
\bigcap_{K\in \mathrm{EP}(M)}\mathrm{An}(K)$ \\
& $=$ & $\mathrm{An}(\sum_{K\in \mathrm{EP}(M)}K)$ & $=$ & $\mathrm{An}(%
\mathrm{EPcorad}(M)).$%
\end{tabular}%
\end{equation*}%
If moreover $M$ is self-cogenerator, then
\begin{equation*}
\mathrm{EPcorad}(M)=\mathrm{KeAn}(\mathrm{EPcorad}(M))=\mathrm{Ke}(\mathrm{%
Prad}(\mathrm{E}_{M}^{\mathcal{C}})).\blacksquare
\end{equation*}
\end{Beweis}

\begin{corollary}
\label{ME-semi}Let $M$ be intrinsically injective self-cogenerator. If $%
\mathrm{E}_{M}^{\mathcal{C}}$ is right Noetherian, then%
\begin{equation*}
M=\mathrm{EPcorad}(M)\Leftrightarrow M\text{ is }\mathrm{E}\text{-semiprime}.
\end{equation*}
\end{corollary}

\begin{Beweis}
Under the assumptions and Proposition \ref{P-rad} we have: $M=\mathrm{EPcorad%
}(M)\Rightarrow \mathrm{Prad}(\mathrm{E}_{M}^{\mathcal{C}})=\mathrm{An}(%
\mathrm{EPcorad}(M))=\mathrm{An}(M)=0,$ i.e. $\mathrm{E}_{M}^{\mathcal{C}}$
is semiprime; on the otherhand $\mathrm{E}_{M}^{\mathcal{C}}$ semiprime $%
\Rightarrow \mathrm{EPcorad}(M)=\mathrm{Ke}(\mathrm{Prad}(\mathrm{E}_{M}^{%
\mathcal{C}}))=\mathrm{Ke}(0)=M.\blacksquare $
\end{Beweis}

\begin{remark}
Let $_{A}\mathcal{C}$ be locally projective and $M$ be right $\mathcal{C}$%
-comodule. A sufficient condition for $\mathrm{E}_{M}^{\mathcal{C}}$ to be
right Noetherian, so that the results of Proposition \ref{P-rad} and
Corollary \ref{ME-semi} follow, is that $M$ is Artinian and self-injective
(see \ref{End} (1)). We recall here also that in case $A_{A}$ is Artinian,
every finitely generated right $\mathcal{C}$-comodule has finite length by
\cite[Corollary 2.25 (4)]{Abu03}.
\end{remark}

\begin{proposition}
\label{th}Let $\theta :L\rightarrow M$ be an isomorphism of right $\mathcal{C%
}$-comodules. Then we have bijections%
\begin{equation}
\mathrm{EP}(L)\leftrightarrow \mathrm{EP}(M),\text{ }\mathrm{ESP}%
(L)\leftrightarrow \mathrm{ESP}(M),\text{ }\mathrm{CEP}(L)\leftrightarrow
\mathrm{CEP}(M),\text{ }\mathrm{CESP}(L)\leftrightarrow \mathrm{CESP}(M).
\label{th-corresp}
\end{equation}%
In particular
\begin{equation}
\theta (\mathrm{EPcorad}(L))=\mathrm{EPcorad}(M).  \label{th-Pcorad}
\end{equation}%
If moreover $L,M$ are self-injective, then we have bijections between the
class of $\mathrm{E}$-prime \emph{(}respectively $\mathrm{E}$-semiprime,
completely $\mathrm{E}$-prime, completely $\mathrm{E}$-semiprime\emph{) }$%
\mathcal{C}$-subcomodules of $L$ and the class of $\mathrm{E}$-prime \emph{(}%
respectively $\mathrm{E}$-semiprime, completely $\mathrm{E}$-prime,
completely $\mathrm{E}$-semiprime\emph{) }$\mathcal{C}$-subcomodules of $M.$
\end{proposition}

\begin{Beweis}
Sine $\theta $ is an isomorphism in $\mathbb{M}^{\mathcal{C}},$ we have an
isomorphism of $R$-algebras
\begin{equation*}
\widetilde{\theta }:\mathrm{E}_{M}^{\mathcal{C}}\rightarrow \mathrm{E}_{L}^{%
\mathcal{C}},\text{ }f\mapsto \lbrack \theta ^{-1}\circ f\circ \theta ].
\end{equation*}%
The result follows then from the fact that for every fully invariant $%
\mathcal{C}$-subcomodule $0\neq H\subseteq L$ (respectively $0\neq
K\subseteq M$), $\widetilde{\theta }$ induces an isomorphism of $R$-algebras
\begin{equation*}
\mathrm{E}_{M}^{\mathcal{C}}/\mathrm{An}(\theta (H))\simeq \mathrm{E}_{L}^{%
\mathcal{C}}/\mathrm{An}(H)\text{ (respectively }\mathrm{E}_{L}^{\mathcal{C}%
}/\mathrm{An}(\theta ^{-1}(K))\simeq \mathrm{E}_{M}^{\mathcal{C}}/\mathrm{An}%
(K)\text{).}\blacksquare
\end{equation*}
\end{Beweis}

\begin{remark}
\label{ev}Let $L$ be a non-zero right $\mathcal{C}$-comodule and $\theta
:L\longrightarrow M$ be a $\mathcal{C}$-colinear map. If $\theta $ is not
bijective, then it is NOT evident that we have the correspondences (\ref%
{th-corresp}).
\end{remark}

\qquad Despite Remark \ref{ev} we have

\begin{proposition}
\label{comod-internal}Let $M$ be self-injective and $0\neq L\subseteq M$ be
a fully invariant non-zero $\mathcal{C}$-subcomodule. Then
\begin{equation*}
\begin{tabular}{lllllll}
$\mathcal{C}_{f.i.}(L)\cap \mathrm{EP}(M)$ & $=$ & $\mathrm{EP}(L)$ & $;$ & $%
\mathcal{C}_{f.i.}(L)\cap \mathrm{CEP}(M)$ & $=$ & $\text{ }\mathrm{CEP}(L)$
\\
$\mathcal{C}_{f.i.}(L)\cap \mathrm{ESP}(M)$ & $=$ & $\mathrm{ESP}(L)$ & $;$
& $\mathcal{C}_{f.i.}(L)\cap \mathrm{CESP}(M)$ & $=$ & $\mathrm{CESP}(L).$%
\end{tabular}%
\end{equation*}
\end{proposition}

\begin{Beweis}
Assume $M$ to be self-injective (so that $L$ is also self-injective). Let $%
0\neq K\subseteq L$ be an arbitrary non-zero fully invariant $\mathcal{C}$%
-subcomodule (so that $K\subseteq M$ is also fully invariant). The result
follows then directly from the definitions and the canonical isomorphisms of
$R$-algebras
\begin{equation*}
\mathrm{E}_{M}^{\mathcal{C}}/\mathrm{An}_{\mathrm{E}_{M}^{\mathcal{C}%
}}(K)\simeq \mathrm{E}_{K}^{\mathcal{C}}\simeq \mathrm{E}_{L}^{\mathcal{C}}/%
\mathrm{An}_{\mathrm{E}_{L}^{\mathcal{C}}}(K).\blacksquare
\end{equation*}
\end{Beweis}

\begin{corollary}
\label{s.i.-internal}Let $M$ be self-injective and $0\neq L\subseteq M$ be a
non-zero fully invariant $\mathcal{C}$-subcomodule. Then $L\in \mathrm{EP}%
(M) $ \emph{(}respectively $L\in \mathrm{ESP}(M),$ $L\in \mathrm{CEP}(M),$ $%
L\in \mathrm{CESP}(M)$\emph{)} if and only if $L$ is $\mathrm{E}$-prime
\emph{(}respectively $\mathrm{E}$-semiprime, completely $\mathrm{E}$-prime,
completely $\mathrm{E}$-semiprime\emph{)}.
\end{corollary}

\subsection*{Sufficient and necessary conditions}

\qquad Given a fully invariant non-zero $\mathcal{C}$-subcomodule $%
K\subseteq M,$ we give sufficient and necessary conditions for $\mathrm{An}%
(K)\vartriangleleft \mathrm{E}_{M}^{\mathcal{C}}$ to be prime (respectively
semiprime, completely prime, completely semiprime). These generalize the
conditions given in \cite{YDZ90} for the dual algebras of a coalgebra over a
base field to be prime (respectively semiprime, domain).

\begin{proposition}
\label{iff-coprime}Let $0\neq K\subseteq M$ be a non-zero fully invariant $%
\mathcal{C}$-subcomodule. A sufficient condition for $K$ to be in $\mathrm{EP%
}(M)$ is that $K=Kf\mathrm{E}_{M}^{\mathcal{C}}$ $\forall f\in \mathrm{E}%
_{M}^{\mathcal{C}}\backslash \mathrm{An}(K),$ where the later is also
necessary in case $M$ is self-cogenerator \emph{(}or $M$ is self-injective
and $K$ is coretractable\emph{)}.
\end{proposition}

\begin{Beweis}
Let $I,J\vartriangleleft \mathrm{E}_{M}^{\mathcal{C}}$ with $IJ\subseteq
\mathrm{An}(K).$ Suppose $I\nsubseteq \mathrm{An}(K)$ and pick some $f\in
I\backslash \mathrm{An}(K).$ By assumption $K=Kf\mathrm{E}_{M}^{\mathcal{C}}$
and it follows then that $KJ=(Kf\mathrm{E}_{M}^{\mathcal{C}})J\subseteq
K(IJ)=0,$ i.e. $J\subseteq \mathrm{An}(K).$

On the otherhand, assume $M$ is self-cogenerator \emph{(}or $M$ is
self-injective and $K$ is coretractable\emph{)}. Suppose there exists some $%
f\in \mathrm{E}_{M}^{\mathcal{C}}\backslash \mathrm{An}(K),$ such that $H:=Kf%
\mathrm{E}_{M}^{\mathcal{C}}\subsetneqq K\neq 0.$ Then obviously $(\mathrm{E}%
_{M}^{\mathcal{C}}f\mathrm{E}_{M}^{\mathcal{C}})\mathrm{An}(H)\subseteq
\mathrm{An}(K),$ whereas our assumptions and Remarks \ref{inj-cor} (3) $\&$
(4) imply that $\mathrm{E}_{M}^{\mathcal{C}}f\mathrm{E}_{M}^{\mathcal{C}%
}\nsubseteq \mathrm{An}(K)$ and $\mathrm{An}(H)\nsubseteq \mathrm{An}(K)$
(i.e. $\mathrm{An}(K)$ is not prime).$\blacksquare $
\end{Beweis}

\begin{proposition}
\label{iff-semi}Let $0\neq K\subseteq M$ be a non-zero fully invariant $%
\mathcal{C}$-subcomodule. A sufficient condition for $K$ to be in $\mathrm{%
ESP}(M)$ is that $Kf=Kf\mathrm{E}_{M}^{\mathcal{C}}f$ $\forall f\in \mathrm{E%
}_{M}^{\mathcal{C}}\backslash \mathrm{An}(K),$ where the later is also
necessary in case $M$ is self-cogenerator.
\end{proposition}

\begin{Beweis}
Let $I^{2}\subseteq \mathrm{An}(K)$ for some $I\vartriangleleft \mathrm{E}%
_{M}^{\mathcal{C}}.$ Suppose $I\nsubseteq \mathrm{An}(K)$ and pick some $%
f\in I\backslash \mathrm{An}(K).$ Then $0\neq Kf\neq Kf\mathrm{E}_{M}^{%
\mathcal{C}}f\subseteq KI^{2}=0,$ a contradiction. So $I\subseteq \mathrm{An}%
(K).$

On the otherhand, assume that $M$ is self-cogenerator. Suppose there exists
some $f\in \mathrm{E}_{M}^{\mathcal{C}}\backslash \mathrm{An}(K)$ with $Kf%
\mathrm{E}_{M}^{\mathcal{C}}f\varsubsetneqq Kf\neq 0.$ By assumptions and
Remark \ref{inj-cor} (3), there exists some $g\in \mathrm{An(}Kf\mathrm{E}%
_{M}^{\mathcal{C}}f)\backslash \mathrm{An}(Kf)$ and it follows then that $J:=%
\mathrm{E}_{M}^{\mathcal{C}}(fg)\mathrm{E}_{M}^{\mathcal{C}}\nsubseteq
\mathrm{An}(K)$ while $J^{2}\subseteq \mathrm{An}(K)$ (i.e. $\mathrm{An}%
(K)\vartriangleleft \mathrm{E}_{M}^{\mathcal{C}}$ is not semiprime).$%
\blacksquare $
\end{Beweis}

\begin{proposition}
\label{iff-c-coprime}Let $0\neq K\subseteq M$ be a non-zero fully invariant $%
\mathcal{C}$-subcomodule. A sufficient condition for $K$ to be in $\mathrm{%
CEP}(M)$ is that $K=Kf$ $\forall f\in \mathrm{E}_{M}^{\mathcal{C}}\backslash
\mathrm{An}(K),$ where the later is also necessary in case $M$ is
self-cogenerator \emph{(}or $M$ is self-injective and $K$ is coretractable%
\emph{)}.
\end{proposition}

\begin{Beweis}
\begin{enumerate}
\item Let $fg\in \mathrm{An}(K)$ for some $f,g\in \mathrm{E}_{M}^{\mathcal{C}%
}$ and suppose $f\notin \mathrm{An}(K).$ The assumption $K=Kf$ implies then
that $Kg=(Kf)g=K(fg)=$ $0,$ i.e. $g\in \mathrm{An}(K).$

On the otherhand, assume $M$ is self-cogenerator \emph{(}or $M$ is
self-injective and $K$ is coretractable\emph{)}. Suppose $Kf\varsubsetneqq
K\neq 0$ for some $f\in \mathrm{E}_{M}^{\mathcal{C}}\backslash \mathrm{An}%
(K) $. By assumptions and Remarks \ref{inj-cor} (3) $\&$ (4) there exists
some $g\in \mathrm{An}(Kf)\backslash \mathrm{An}(K)$ with $fg\in \mathrm{An}%
(K)\ $(i.e. $\mathrm{An}(K)\vartriangleleft \mathrm{E}_{M}^{\mathcal{C}}$ is
not completely prime).$\blacksquare $
\end{enumerate}
\end{Beweis}

\begin{proposition}
\label{iff-c-semi}Let $0\neq K\subseteq M$ be a non-zero fully invariant $%
\mathcal{C}$-subcomodule. A sufficient condition for $K$ to be in $\mathrm{%
CESP}(M)$ is that $Kf=Kf^{2}$ for every $f\in \mathrm{E}_{M}^{\mathcal{C}%
}\backslash \mathrm{An}(K),$ where the later is also necessary in case $M$
is self-cogenerator.
\end{proposition}

\begin{Beweis}
Let $f\in \mathrm{E}_{M}^{\mathcal{C}}$ be such that $f^{2}\in \mathrm{An}%
(K).$ The assumption $K=Kf$ implies then that $Kf=Kf^{2}=0,$ i.e. $f\in
\mathrm{An}(K).$ On the otherhand, assume $M$ is self-cogenerator. Suppose
that $Kf^{2}\varsubsetneqq Kf\neq 0$ for some $f\in \mathrm{E}_{M}^{\mathcal{%
C}}\backslash \mathrm{An}(K).$ By assumptions and Remark \ref{inj-cor} (3),
there exists some $g\in \mathrm{An}(Kf^{2})\backslash \mathrm{An}(Kf).$ Set
\begin{equation*}
h:=%
\begin{cases}
fgf, & \text{in case }fgf\notin \mathrm{An}(K); \\
fg, & \text{otherwise.}%
\end{cases}%
\end{equation*}%
So $h^{2}\in \mathrm{An}(K)$ while $h\notin \mathrm{An}(K)$ (i.e. $\mathrm{An%
}(K)\vartriangleleft \mathrm{E}_{M}^{\mathcal{C}}$ is not completely
semiprime).$\blacksquare $
\end{Beweis}

\qquad The proof of the following result can be obtained directly from the
proofs of the previous four propositions by replacing $K$ with $M.$

\begin{theorem}
\label{iff-M}

\begin{enumerate}
\item $M$ is \emph{(}completely\emph{)} $\mathrm{E}$-prime, if $M=Mf\mathrm{E%
}_{M}^{\mathcal{C}}$ \emph{(}$M=Mf$\emph{)} for every $0\neq f\in \mathrm{E}%
_{M}^{\mathcal{C}}.$ If $M$ is coretractable, then $M$ is \emph{(}completely%
\emph{)} $\mathrm{E}$-prime if and only if $M=Mf\mathrm{E}_{M}^{\mathcal{C}}$
\emph{(}$M=Mf$\emph{)} for every $0\neq f\in \mathrm{E}_{M}^{\mathcal{C}}.$

\item $M$ is \emph{(}completely\emph{)} $\mathrm{E}$-semiprime,\emph{\ }if $%
Mf=Mf\mathrm{E}_{M}^{\mathcal{C}}f$ \emph{(}$Mf=Mf^{2}$\emph{)} for every $%
0\neq f\in \mathrm{E}_{M}^{\mathcal{C}}.$ If $M$ is self-cogenerator, then $%
M $ is \emph{(}completely\emph{)} $\mathrm{E}$-semiprime\emph{\ }if and only
if $Mf=Mf\mathrm{E}_{M}^{\mathcal{C}}f$ \emph{(}$Mf=Mf^{2}$\emph{)} for
every $0\neq f\in \mathrm{E}_{M}^{\mathcal{C}}.$
\end{enumerate}
\end{theorem}

\subsection*{\textrm{E}-Prime versus simple}

\qquad In what follows we show that \textrm{E}-prime comodules generalize
simple comodules.

\begin{theorem}
\label{r-simple-E}A sufficient condition for $\mathrm{E}_{M}^{\mathcal{C}}$
to be right simple \emph{(}a division ring\emph{)} is that $M$ is simple,
where the later is also necessary in case $M$ is self-cogenerator.
\end{theorem}

\begin{Beweis}
If $M$ is simple, then $\mathrm{E}_{M}^{\mathcal{C}}:=\mathrm{End}^{\mathcal{%
C}}(M)^{op}$ is a Division ring by Schur's Lemma.

On the otherhand, assume $M$ to be self-cogenerator. Let $K\subseteq M$ be a
$\mathcal{C}$-subcomodule and consider the right ideal $\mathrm{An}%
(K)\vartriangleleft _{r}\mathrm{E}_{M}^{\mathcal{C}}.$ If $\mathrm{E}_{M}^{%
\mathcal{C}}$ is right simple, then $\mathrm{An}(K)=(0_{\mathrm{E}_{M}^{%
\mathcal{C}}})$ so that $K=\mathrm{KeAn}(K)=\mathrm{Ke}(0_{\mathrm{E}_{M}^{%
\mathcal{C}}})=M;$ or $\mathrm{An}(K)=\mathrm{E}_{M}^{\mathcal{C}}$ so that $%
K=\mathrm{KeAn}(K)=\mathrm{Ke}(\mathrm{E}_{M}^{\mathcal{C}})=(0_{M}).$
Consequently $M$ is simple.$\blacksquare $
\end{Beweis}

\begin{theorem}
\label{E-simple}A sufficient condition for $\mathrm{E}_{M}^{\mathcal{C}}$ to
be simple, in case $M$ is intrinsically injective and $\mathrm{E}_{M}^{%
\mathcal{C}}$ is right Noetherian, is to assume that $M$ has no non-trivial
fully invariant $\mathcal{C}$-subcomodules, where the later is also
necessary if $M$ is self-cogenerator.
\end{theorem}

\begin{Beweis}
The proof is similar to that of Theorem \ref{r-simple-E} replacing right
ideals of $\mathrm{E}_{M}^{\mathcal{C}}$ by two-sided ideals and arbitrary $%
\mathcal{C}$-subcomodules of $M$ with fully invariant ones.$\blacksquare $
\end{Beweis}

\begin{notation}
Consider the non-zero right $\mathcal{C}$-comodule $M.$ With $\mathcal{S}(M)$
($\mathcal{S}_{f.i.}(M)$) we denote the class of simple $\mathcal{C}$%
-subcomodules of $M$ (fully invariant $\mathcal{C}$-subcomodules of $M$ with
no non-trivial fully invariant $\mathcal{C}$-subcomodules).
\end{notation}

\begin{corollary}
\label{max-simple}Let $M$ be self-injective self-cogenerator and $0\neq
K\subseteq M$ be a fully invariant non-zero $\mathcal{C}$-subcomodule. Then

\begin{enumerate}
\item $K\in \mathcal{S}(M)\Leftrightarrow \mathrm{An}(K)\in \mathrm{Max}_{r}(%
\mathrm{E}_{M}^{\mathcal{C}});$

\item If $\mathrm{E}_{M}^{\mathcal{C}}$ is right Noetherian, then $\mathcal{S%
}_{f.i.}(M)\Leftrightarrow \mathrm{An}(K)\in \mathrm{Max}(\mathrm{E}_{M}^{%
\mathcal{C}}).$
\end{enumerate}
\end{corollary}

\begin{Beweis}
Recall that, since $M$ is self-injective self-cogenerator and $K\subseteq M$
is fully invariant, $K$ is also self-injective self-cogenerator. The result
follows then from Theorems \ref{r-simple-E} and \ref{E-simple} applied to
the $R$-algebra $\mathrm{E}_{K}^{\mathcal{C}}\simeq \mathrm{E}_{M}^{\mathcal{%
C}}/\mathrm{An}(K).\blacksquare $
\end{Beweis}

\begin{lemma}
\label{simple-max}Let $M$ be intrinsically injective self-cogenerator and
assume $\mathrm{E}_{M}^{\mathcal{C}}$ to be right Noetherian. Then the order
reversing mappings (\ref{An-map-Ke}) give a bijection%
\begin{equation}
\mathcal{S}(M)\longleftrightarrow \mathrm{Max}_{r}(\mathrm{E}_{M}^{\mathcal{C%
}})\text{ and }\mathcal{S}_{f.i.}(M)\longleftrightarrow \mathrm{Max}(\mathrm{%
E}_{M}^{\mathcal{C}}).  \label{simple-An-Max}
\end{equation}
\end{lemma}

\begin{Beweis}
Let $K\in \mathcal{S}(M)$ ($K\in \mathcal{S}_{f.i.}(M)$) and consider the
proper right ideal $\mathrm{An}(K)\subsetneqq \mathrm{E}_{M}^{\mathcal{C}}.$
If $\mathrm{An}(K)\subseteq I\subseteq \mathrm{E}_{M}^{\mathcal{C}},$ for
some right (two-sided) ideal $I\subseteq \mathrm{E}_{M}^{\mathcal{C}},$ then
$\mathrm{Ke}(I)\subseteq \mathrm{KeAn}(K)=K$ and it follows from the
assumption $K\in \mathcal{S}(M)$ ($K\in \mathcal{S}_{f.i.}(M)$) that $%
\mathrm{Ke}(I)=0$ so that $I=\mathrm{AnKe}(I)=\mathrm{E}_{M}^{\mathcal{C}};$
or $\mathrm{Ke}(I)=K$ so that $I=\mathrm{AnKe}(I)=\mathrm{An}(K).$ This
means that $\mathrm{An}(K)\in \mathrm{Max}_{r}(\mathrm{E}_{M}^{\mathcal{C}})$
($\mathrm{An}(K)\in \mathrm{Max}(\mathrm{E}_{M}^{\mathcal{C}})$).

On the otherhand, let $I\in \mathrm{Max}_{r}(\mathrm{E}_{M}^{\mathcal{C}})$ (%
$I\in \mathrm{Max}(\mathrm{E}_{M}^{\mathcal{C}})$) and consider the non-zero
$\mathcal{C}$-subcomodule $0\neq \mathrm{Ke}(I)\subseteq M.$ If $K\subseteq
\mathrm{Ke}(I)$ for some (fully invariant) $\mathcal{C}$-subcomodule $%
K\subseteq M,$ then $I\subseteq \mathrm{AnKe}(I)\subseteq \mathrm{An}%
(K)\subseteq \mathrm{E}_{M}^{\mathcal{C}}$ and it follows by the maximality
of $I$ that $\mathrm{An}(K)=\mathrm{E}_{M}^{\mathcal{C}}$ so that $K=\mathrm{%
KeAn}(K)=0;$ or $\mathrm{An}(K)=I$ so that $K=\mathrm{KeAn}(K)=\mathrm{Ke}%
(I).$ Consequently $\mathrm{Ke}(I)\in \mathcal{S}(M)$ ($K\in \mathcal{S}%
_{f.i.}(M)$). Since $M$ is intrinsically injective self-cogenerator, $%
\mathrm{Ke}(-)$ and $\mathrm{An}(-)$ are injective by \ref{An-Ke} and we are
done.$\blacksquare $
\end{Beweis}

\begin{lemma}
\label{essential}Let $A$ be left perfect and $_{A}\mathcal{C}$ be locally
projective.

\begin{enumerate}
\item The non-zero right $\mathcal{C}$-comodule contains a simple $\mathcal{C%
}$-subcomodule.

\item $\mathrm{Soc}(M)\vartriangleleft _{e}M$ \emph{(}an essential $\mathcal{%
C}$-subcomodule\emph{)}.
\end{enumerate}
\end{lemma}

\begin{Beweis}
Let $_{A}A$ be perfect and $_{A}\mathcal{C}$ be locally projective.

\begin{enumerate}
\item By \cite[Corollary 2.25]{Abu03} $M$ satisfies the descending chain
condition on finitely generated non-zero $\mathcal{C}$-subcomodules, which
turn out to be finitely generated right $A$-modules, hence $M$ contains a
non-zero simple $\mathcal{C}$-subcomodule.

\item Let $M$ be a non-zero right $\mathcal{C}$-comodule. For every $%
\mathcal{C}$-subcomodule $0\neq K\subseteq M$ we have $K\cap \mathrm{Soc}(M)=%
\mathrm{Soc}(K)\neq 0,$ by (1).$\blacksquare $
\end{enumerate}
\end{Beweis}

\begin{proposition}
\label{Jac}We have%
\begin{equation}
\mathrm{Jac}(\mathrm{E}_{M}^{\mathcal{C}})=\mathrm{An}(\mathrm{Soc}(M))\text{
and }\mathrm{Soc}(M)=\mathrm{Ke}(\mathrm{Jac}(\mathrm{E}_{M}^{\mathcal{C}})),
\end{equation}%
if any of the following conditions holds:

\begin{enumerate}
\item $M$ is intrinsically injective self-cogenerator with $\mathrm{E}_{M}^{%
\mathcal{C}}$ right Noetherian;

\item $_{A}\mathcal{C}$ is locally projective and $M$ is Artinian
self-injective self cogenerator;

\item $A$ is left perfect, $_{A}\mathcal{C}$ is locally projective and $M$
is self-injective self-cogenerator\emph{.}
\end{enumerate}
\end{proposition}

\begin{Beweis}
\begin{enumerate}
\item By Lemma \ref{simple-max} we have%
\begin{equation*}
\begin{tabular}{lll}
$\mathrm{Jac}(\mathrm{E}_{M}^{\mathcal{C}})$ & $=$ & $\bigcap \{Q\mid
Q\vartriangleleft _{r}\mathrm{E}_{M}^{\mathcal{C}}$ is a maximal right ideal$%
\}$ \\
& $=$ & $\bigcap \{\mathrm{AnKe}(Q)\mid Q\vartriangleleft _{r}\mathrm{E}%
_{M}^{\mathcal{C}}$ is a maximal right ideal$\}$ \\
& $=$ & $\bigcap \{\mathrm{AnKe}(\mathrm{An}(K))\mid K\subseteq M\text{ is a
simple }\mathcal{C}\text{-subcomodule}\}$ \\
& $=$ & $\bigcap \{\mathrm{An}(K)\mid K\subseteq M\text{ is a simple }%
\mathcal{C}\text{-subcomodule}\}$ \\
& $=$ & $\mathrm{An}(\sum \{K\mid K\subseteq M\text{ is a simple }\mathcal{C}%
\text{-subcomodule}\})$ \\
& $=$ & $\mathrm{An}(\mathrm{Soc}(M)).$%
\end{tabular}%
\end{equation*}%
Since $M$ is self-cogenerator, we have $\mathrm{Soc}(M)=\mathrm{KeAn}(%
\mathrm{Soc}(M))=\mathrm{Ke}(\mathrm{Jac}(\mathrm{E}_{M}^{\mathcal{C}})).$

\item Since $M$ is Artinian and self-injective in $\mathbb{M}^{\mathcal{C}%
}=\sigma \lbrack _{^{\ast }\mathcal{C}}\mathcal{C}],$ we conclude that $%
\mathrm{E}_{M}^{\mathcal{C}}:=\mathrm{End}^{\mathcal{C}}(M)^{op}=\mathrm{End}%
(_{^{\ast }\mathcal{C}}M)$ is right Noetherian by Proposition \ref{End} (2).
The result follows then by (1).

\item Since $A$ is left perfect and $_{A}\mathcal{C}$ is locally projective,
$\mathrm{Soc}(M)\vartriangleleft _{e}M$ is an essential $\mathcal{C}$%
-subcomodule by Lemma \ref{essential} (2) and it follows then, since $M$ is
self-injective, that%
\begin{equation*}
\begin{tabular}{llll}
$\mathrm{Jac}(\mathrm{E}_{M}^{\mathcal{C}})$ & $=\mathrm{Jac}(\mathrm{End}%
(_{^{\ast }\mathcal{C}}M)^{op})$ & $=\mathrm{Hom}_{^{\ast }\mathcal{C}}(M/%
\mathrm{Soc}(M),M)$ & (\cite[22.1 (5)]{Wis91}) \\
& $=\mathrm{Hom}^{\mathcal{C}}(M/\mathrm{Soc}(M),M)$ & $\simeq \mathrm{An}(%
\mathrm{Soc}(M)).$ &
\end{tabular}%
\end{equation*}%
Since $M$ is self-cogenerator, we have moreover%
\begin{equation*}
\mathrm{Soc}(M)=\mathrm{KeAn}(\mathrm{Soc}(M))=\mathrm{Ke}(\mathrm{Jac}(%
\mathrm{E}_{M}^{\mathcal{C}})).\blacksquare
\end{equation*}
\end{enumerate}
\end{Beweis}

\begin{corollary}
\label{semisimple-semiprimitive}If any of the three conditions in
Proposition \ref{Jac} holds, then we have%
\begin{equation*}
M\text{ is semisimple}\Leftrightarrow \mathrm{E}_{M}^{\mathcal{C}}\text{ is
semiprimitive.}
\end{equation*}
\end{corollary}

\begin{Beweis}
By assumptions and Proposition \ref{Jac} we have $\mathrm{Jac}(\mathrm{E}%
_{M}^{\mathcal{C}})=\mathrm{An}(\mathrm{Soc}(M))$ and $\mathrm{Soc}(M)=%
\mathrm{Ke}(\mathrm{Jac}(\mathrm{E}_{M}^{\mathcal{C}})).$ Hence, $M$
semisimple $\Rightarrow \mathrm{Jac}(\mathrm{E}_{M}^{\mathcal{C}})=\mathrm{An%
}(\mathrm{Soc}(M))=\mathrm{An}(M)=0,$ i.e. $\mathrm{E}_{M}^{\mathcal{C}}$ is
semiprimitive; on the otherhand $\mathrm{E}_{M}^{\mathcal{C}}$ semiprimitive
implies $\mathrm{Soc}(M)=\mathrm{Ke}(\mathrm{Jac}(\mathrm{E}_{M}^{\mathcal{C}%
}))=\mathrm{Ke}(0)=M,$ i.e. $M$ is semisimple.$\blacksquare $
\end{Beweis}

\subsection*{\textrm{E}-Prime versus irreducible}

\qquad In what follows we clarify the relation between $\mathrm{E}$-prime
and irreducible comodules.

\begin{remark}
\label{int-sem}Let $\{K_{\lambda }\}_{\Lambda }$ be a family of non-zero
fully invariant $\mathcal{C}$-subcomodules of $M$ and consider the fully
invariant $\mathcal{C}$-subcomodule $K:=\sum\limits_{\lambda \in \Lambda
}K_{\lambda }\subseteq M.$ If $K_{\lambda }\in \mathrm{EP}(M)$ $(K_{\lambda
}\in \mathrm{CEP}(M))$ for every $\lambda \in \Lambda ,$ then $\mathrm{An}%
(K)=\bigcap\limits_{\lambda \in \Lambda }\mathrm{An}(K_{\lambda })\ $is an
intersection of (completely) prime ideals, hence a (completely) semiprime
ideal, i.e. $K\in \mathrm{ESP}(M)$ ($K\in \mathrm{CESP}(M)$). If $M$ is
self-injective, then we conclude that an arbitrary sum of (completely) $%
\mathrm{E}$-prime $\mathcal{C}$-subcomodules of $M$ is in general
(completely) $\mathrm{E}$-semiprime.
\end{remark}

\qquad Despite Remark \ref{int-sem} we have the following result (which is
most interesting in case $K=M$):

\begin{proposition}
\label{family}Let $\{K_{\lambda }\}_{\Lambda }$ be a family of non-zero
fully invariant $\mathcal{C}$-subcomodules of $M,$ such that for any $\gamma
,\delta \in \Lambda $ either $K_{\gamma }$ $\subseteq K_{\delta }$ or $%
K_{\delta }\subseteq K_{\gamma },$ and consider the fully invariant $%
\mathcal{C}$-subcomodule $K:=\sum\limits_{\lambda \in \Lambda }K_{\lambda
}=\bigcup\limits_{\lambda \in \Lambda }K_{\lambda }\subseteq M.$ If $%
K_{\lambda }\in \mathrm{EP}(M)$ \emph{(}$K_{\lambda }\in \mathrm{CEP}(M)$%
\emph{)} for every $\lambda \in \Lambda ,$ then $K\in \mathrm{EP}(M)$ \emph{(%
}$K\in \mathrm{CEP}(M)$\emph{).}

\begin{Beweis}
Let $I,J\vartriangleleft \mathrm{E}_{M}^{\mathcal{C}}$ be such that $%
IJ\subseteq \mathrm{An}(K)=\bigcap\limits_{\lambda \in \Lambda }\mathrm{An}%
(K_{\lambda })$ and suppose $I\nsubseteq \mathrm{An}(K).$ Pick some $\lambda
_{0}\in \Lambda $ with $I\nsubseteq \mathrm{An}(K_{\lambda _{0}}).$ By
assumption $\mathrm{An}(K_{\lambda _{0}})\vartriangleleft \mathrm{E}_{M}^{%
\mathcal{C}}$ is prime and $IJ\subseteq \mathrm{An}(K_{\lambda _{0}}),$ so $%
J\subseteq \mathrm{An}(K_{\lambda _{0}}).$ We \textbf{claim} that $%
J\subseteq \bigcap\limits_{\lambda \in \Lambda }\mathrm{An}(K_{\lambda }):$
Let $\lambda \in \Lambda $ be arbitrary. If $K_{\lambda }\subseteq
K_{\lambda _{0}},$ then $J\subseteq \mathrm{An}(K_{\lambda _{0}})\subseteq
\mathrm{An}(K_{\lambda }).$ On the other hand, if $K_{\lambda _{0}}\subseteq
K_{\lambda }$ and $J\nsubseteq \mathrm{An}(K_{\lambda }),$ then the
primeness of $\mathrm{An}(K_{\lambda })$ implies that $I\subseteq \mathrm{An}%
(K_{\lambda })\subseteq \mathrm{An}(K_{\lambda _{0}}),$ a contradiction. So $%
J\subseteq \bigcap\limits_{\lambda \in \Lambda }\mathrm{An}(K_{\lambda })=%
\mathrm{An}(K).$ Consequently $\mathrm{An}(K)\vartriangleleft \mathrm{E}%
_{M}^{\mathcal{C}}$ is a prime ideal, i.e. $K\in \mathrm{EP}(M).$

In case $\mathrm{An}(K_{\lambda })\vartriangleleft \mathrm{E}_{M}^{\mathcal{C%
}}$ is completely prime for every $\lambda \in \Lambda ,$ then replacing
ideals in the argument above with elements yields that $\mathrm{An}%
(K)\vartriangleleft \mathrm{E}_{M}^{\mathcal{C}}$ is a completely prime
ideal, i.e. $K\in \mathrm{CEP}(M).\blacksquare $
\end{Beweis}
\end{proposition}

\begin{remark}
If $M$ is self-injective and the subcomodule $K_{\lambda }$ in Proposition %
\ref{family} are \emph{(}completely\emph{)} $\mathrm{E}$-prime, then $K$ is
\emph{(}completely\emph{)} $\mathrm{E}$-prime (recall that we have in this
case an isomorphism of algebras $\mathrm{E}_{M}^{\mathcal{C}}/\mathrm{Ann}%
(K_{\lambda })\simeq \mathrm{E}_{K_{\lambda }}^{\mathcal{C}}$).
\end{remark}

\begin{proposition}
\label{full-inv-dec}Let $M$ be self-cogenerator and $K\in \mathrm{EP}(M).$
Then $K$ admits no decomposition as an internal direct sum of non-trivial
fully invariant $\mathcal{C}$-subcomodules.
\end{proposition}

\begin{Beweis}
Let $K\subseteq M$ be a fully invariant $\mathcal{C}$-subcomodule with $%
\mathrm{An}(K)\vartriangleleft \mathrm{E}_{M}^{\mathcal{C}}$ a prime ideal
and suppose $K=K_{\lambda _{0}}\oplus \sum\limits_{\lambda \neq \lambda
_{0}}K_{\lambda }$ to be a decomposition of $K$ as an internal direct sum of
non-trivial fully invariant $\mathcal{C}$-subcomodules. Consider the
two-sided ideals $I:=\mathrm{An}(K_{\lambda _{0}}),$ $J:=\mathrm{An}%
(\sum\limits_{\lambda \neq \lambda _{0}}K_{\lambda })$ of $\mathrm{E}_{M}^{%
\mathcal{C}},$ so that $IJ\subseteq \mathrm{An}(K).$ If $J\subseteq \mathrm{%
An}(K),$ then $K_{\lambda _{0}}\subseteq K=\mathrm{KeAn}(K)\subseteq \mathrm{%
Ke}(J)=\sum\limits_{\lambda \neq \lambda _{0}}K_{\lambda }$ (a
contradiction). Since $\mathrm{An}(K)\vartriangleleft \mathrm{E}_{M}^{%
\mathcal{C}}$ is prime, $I\subseteq \mathrm{An}(K)$ and we conclude that $K=%
\mathrm{KeAn}(K)\subseteq \mathrm{Ke}(I)=\mathrm{KeAn}(K_{\lambda
_{0}})=K_{\lambda _{0}}$ (a contradiction).$\blacksquare $
\end{Beweis}

\begin{definition}
We call the non-zero right $\mathcal{C}$-comodule $M$ \emph{irreducible},
iff $M$ has a \emph{unique }simple $\mathcal{C}$-subcomodule that is
contained in every $\mathcal{C}$-subcomodule of $M$ (equivalently, iff the
intersection of all non-zero $\mathcal{C}$-subcomodules of $M$ is again
non-zero).
\end{definition}

The following result clarifies, under suitable conditions, the relation
between \textrm{E-}prime and irreducible comodules.

\begin{theorem}
\label{E-coprime-irr}Assume $_{A}\mathcal{C}$ to be locally projective, $M$
to be self-injective self-cogenerator and $\mathrm{End}^{\mathcal{C}}(M)$ to
be commutative. If $M$ is $\mathrm{E}$-prime, then $M$ is irreducible.
\end{theorem}

\begin{Beweis}
If $\mathrm{End}^{\mathcal{C}}(M=\mathrm{End}(_{^{\ast }\mathcal{C}}M)$ is
commutative, then under the assumptions on $M,$ \cite[48.16]{Wis91} yield
that $M$ is a direct sum of irreducible fully invariant $\mathcal{C}$%
-subcomodules. The results follows then by Proposition \ref{full-inv-dec}.$%
\blacksquare $
\end{Beweis}

\section{Fully Coprime (fully cosemiprime) comodules}

\qquad As before, $\mathcal{C}$ is a non-zero $A$-coring with $_{A}\mathcal{C%
}$ flat, $M$ is a non-zero right $\mathcal{C}$-comodule and $\mathrm{E}_{M}^{%
\mathcal{C}}:=\mathrm{End}^{\mathcal{C}}(M)^{op}$ is the ring of $\mathcal{C}
$-colinear endomorphisms of $M$ with the opposite composition of maps.

\begin{punto}
For $R$-submodules $X,Y\subseteq M,$ set%
\begin{equation*}
(X:_{M}^{\mathcal{C}}Y):=\bigcap \{f^{-1}(Y)|\text{ }f\in \mathrm{End}^{%
\mathcal{C}}(M)\text{ and }f(X)=0\}.
\end{equation*}%
If $Y\subseteq M$ is a right $\mathcal{C}$-subcomodule, then $%
f^{-1}(Y)\subseteq M$ is a $\mathcal{C}$-subcomodule for each $f\in \mathrm{E%
}_{M}^{\mathcal{C}},$ being the kernel of the $\mathcal{C}$-colinear map $%
\pi _{Y}\circ f:M\longrightarrow M/Y,$ and it follows then that $(X:_{M}^{%
\mathcal{C}}Y)\subseteq M$ is a right $\mathcal{C}$-subcomodule, being the
intersection of right $\mathcal{C}$-subcomodules of $M.$ If $X\subseteq M$
is fully invariant, i.e. $f(X)\subseteq X$ for every $f\in \mathrm{E}_{M}^{%
\mathcal{C}},$ then $(X:_{M}^{\mathcal{C}}Y)\subseteq M$ is clearly fully
invariant. If $X,Y\subseteq M$ are right $\mathcal{C}$-subcomodules, then
the right $\mathcal{C}$-subcomodule $(X:_{M}^{\mathcal{C}}Y)$ is called the
\emph{internal coproduct }of $X$ and $Y$ in the category $\mathbb{M}^{%
\mathcal{C}}$ of right $\mathcal{C}$-comodules. If $\mathcal{C}_{A}$ is
flat, then the internal coproduct of $\mathcal{C}$-subcomodules of left $%
\mathcal{C}$-comodules can be defined analogously.
\end{punto}

\begin{remark}
The \emph{internal coproduct }of submodules of a given module over a ring
was first introduced by Bican et. al. \cite{BJKN80} to present the notion of
\emph{coprime modules. }The definition was modified in \cite{RRW05}, where
arbitrary submodules are replaced by the fully invariant ones. To avoid any
possible confusion, we refer to coprime modules in the sense of \cite{RRW05}
as \emph{fully coprime modules} and transfer that terminology to \emph{fully
coprime comodules}.
\end{remark}

\begin{definition}
A fully invariant $\mathcal{C}$-subcomodule $0\neq K\subseteq M$ will be
called

\emph{fully }$M$\emph{-coprime, }iff for any two fully invariant $\mathcal{C}
$-subcomodules $X,$ $Y\subseteq M$ with $K\subseteq (X:_{M}^{\mathcal{C}}Y),$
we have $K\subseteq X$ or $K\subseteq Y;$

fully $M$-cosemiprime,\emph{\ }iff for any fully invariant $\mathcal{C}$%
-subcomodule $X\subseteq M$ with $K\subseteq (X:_{M}^{\mathcal{C}}X),$ we
have $K\subseteq X.$

We call $M$ \emph{fully coprime} (\emph{fully cosemiprime}), iff $M$ is
fully $M$-coprime (fully $M$-cosemiprime).
\end{definition}

\subsection*{The fully coprime coradical}

\begin{definition}
We define the \emph{fully coprime spectrum} of $M$\emph{\ }as%
\begin{equation*}
\mathrm{CPSpec}(M):=\{K\mid 0\neq K\subseteq M\text{ is an fully }M\text{%
-coprime }\mathcal{C}\text{-subcomodule}\}.
\end{equation*}%
We define the \emph{fully coprime coradical }of $M$ as%
\begin{equation*}
\mathrm{CPcorad}(M)=\sum_{K\in \mathrm{CPSpec}(M)}K.
\end{equation*}%
Moreover, we set%
\begin{equation*}
\mathrm{CSP}(M):=\{K\mid 0\neq K\subseteq M\text{ is an fully }M\text{%
-cosemiprime }\mathcal{C}\text{-subcomodule}\}.
\end{equation*}
\end{definition}

The fully coprime spectra (fully coprime coradicals) of comodules are
invariant under isomorphisms of comodules:

\begin{proposition}
\label{th-coprime}Let $\theta :L\rightarrow M$ be an isomorphism of $%
\mathcal{C}$-comodules. Then we have bijections
\begin{equation*}
\mathrm{CPSpec}(L)\longleftrightarrow \mathrm{CPSpec}(M)\text{ and }\mathrm{%
CSP}(L)\longleftrightarrow \mathrm{CSP}(M).
\end{equation*}%
In particular
\begin{equation}
\theta (\mathrm{CPcorad}(L))=\mathrm{CPcorad}(M).  \label{th-CPcorad}
\end{equation}
\end{proposition}

\begin{Beweis}
Let $\theta :L\rightarrow M$ be an isomorphism of right $\mathcal{C}$%
-comodules. Let $0\neq H\subseteq L$ be a fully invariant $\mathcal{C}$%
-subcomodule that is fully $L$-coprime and consider the fully invariant $%
\mathcal{C}$-subcomodule $0\neq \theta (H)\subseteq M.$ Let $X,Y\subseteq M$
be two fully invariant $\mathcal{C}$-subcomodules with $\theta (H)\subseteq
(X:_{M}^{\mathcal{C}}Y).$ Then $\theta ^{-1}(X),$ $\theta ^{-1}(Y)\subseteq
L $ are two fully invariant $\mathcal{C}$-subcomodules and $H\subseteq
(\theta ^{-1}(X):_{L}^{\mathcal{C}}\theta ^{-1}(Y)).$ By assumption $H$ is
fully $L$-coprime and we conclude that $H\subseteq \theta ^{-1}(X)$ so that $%
\theta (H)\subseteq X;$ or $H\subseteq \theta ^{-1}(Y)$ so that $\theta
(H)\subseteq Y.$ Consequently $\theta (H)$ is fully $M$-coprime. Analogously
one can show that for any fully invariant fully $M$-coprime $\mathcal{C}$%
-subcomodule $0\neq K\subseteq M,$ the fully invariant $\mathcal{C}$%
-subcomodule $0\neq \theta ^{-1}(K)\subseteq L$ is fully $L$-coprime.

Repeating the proof above with $Y=X,$ one can prove that for any fully $L$%
-cosemiprime (fully $M$-cosemiprime) fully invariant $\mathcal{C}$%
-subcomodule $0\neq H\subseteq L$ (resp. $0\neq K\subseteq M$), the fully
invariant $\mathcal{C}$-subcomodule $0\neq \theta (H)\subseteq M$ ($0\neq
\theta ^{-1}(K)\subseteq L$) is fully $M$-cosemiprime (fully $L$%
-cosemiprime).$\blacksquare $
\end{Beweis}

\begin{remark}
\label{cop-th-not}Let $L$ be a non-zero right $\mathcal{C}$-comodules and $%
\theta :L\rightarrow M$ be a $\mathcal{C}$-colinear map. If $\theta $ is not
bijective, then it is NOT evident that for $K\in \mathrm{CPSpec}(L)$
(respectively $K\in \mathrm{CSP}(L)$) we have $\theta (K)\subseteq \mathrm{%
CPSpec}(M)$ (respectively $\theta (K)\in \mathrm{CSP}(M)$).
\end{remark}

\qquad Despite Remark \ref{cop-th-not} we have

\begin{proposition}
\label{ro-inn}Let $0\neq L\subseteq M$ be a non-zero fully invariant $%
\mathcal{C}$-subcomodule. Then we have%
\begin{equation}
\mathcal{M}_{f.i.}(L)\cap \mathrm{CPSpec}(M)\subseteq \mathrm{CPSpec}(L)%
\text{ and }\mathcal{M}_{f.i.}(L)\cap \mathrm{CSP}(M)\subseteq \mathrm{CSP}%
(L),  \label{it<L}
\end{equation}%
with equality in case $M$ is self injective.
\end{proposition}

\begin{Beweis}
Let $0\neq H\subseteq L$ be a fully invariant $\mathcal{C}$-subcomodule and
assume $H$ to be fully $M$-coprime (fully $M$-cosemiprime). Suppose $%
H\subseteq (X:_{L}^{\mathcal{C}}Y)$ for two (equal) fully invariant $%
\mathcal{C}$-subcomodules $X,$ $Y\subseteq L.$ Since $L\subseteq M$ is a
fully invariant $\mathcal{C}$-subcomodule, it follows that $X,Y$ are also
fully invariant $\mathcal{C}$-subcomodules of $M$ and moreover $(X:_{L}^{%
\mathcal{C}}Y)\subseteq (X:_{M}^{\mathcal{C}}Y).$ By assumption $H$ is fully
$M$-coprime (fully $M$-cosemiprime), and so the inclusions $H\subseteq
(X:_{L}^{\mathcal{C}}Y)\subseteq (X:_{M}^{\mathcal{C}}Y)$ imply $H\subseteq
X $ or $H\subseteq Y.$ Consequently $H$ is fully $L$-coprime (fully $L$%
-cosemiprime). Hence the inclusions in (\ref{it<L}) hold.

Assume now that $M$ is self-injective. Let $0\neq H\subseteq L$ to be an
fully $L$-coprime (fully $L$-cosemiprime) $\mathcal{C}$-subcomodule. Suppose
$X,Y\subseteq M$ are two (equal) fully invariant $\mathcal{C}$-subcomodules
with $H\subseteq (X:_{M}^{\mathcal{C}}Y)$ and consider the fully invariant $%
\mathcal{C}$-subcomodules $X\cap L,$ $Y\cap L\subseteq L.$ Since $M$ is
self-injective, the embedding $\iota :L/X\cap L\hookrightarrow M/X$ induces
a surjective set map
\begin{equation*}
\Phi :\mathrm{Hom}^{\mathcal{C}}(M/X,M)\rightarrow \mathrm{Hom}^{\mathcal{C}%
}(L/X\cap L,M),\text{ }f\mapsto f_{|_{L/X\cap L}}.
\end{equation*}%
Since $L\subseteq M$ is fully invariant, $\Phi $ induces a surjective set
map
\begin{equation}
\Psi :\mathrm{An}_{\mathrm{E}_{M}^{\mathcal{C}}}(X)\rightarrow \mathrm{An}_{%
\mathrm{E}_{L}^{\mathcal{C}}}(X\cap L),\text{ }g\mapsto g_{\mid _{L}},
\label{tel-Phi}
\end{equation}%
which implies that $H\subseteq (X\cap L:_{L}^{\mathcal{C}}Y\cap L).$ By
assumption $H$ is fully $L$-coprime, hence $H\subseteq X\cap L$ so that $%
H\subseteq X;$ or $H\subseteq Y\cap L$ so that $H\subseteq Y.$ Hence $H$ is
fully $M$-coprime (fully $M$-cosemiprime). Consequently the inclusions in (%
\ref{it<L}) become equality.$\blacksquare $
\end{Beweis}

\begin{remark}
\label{ro-internal}Let $0\neq L\subseteq M$ be a non-zero fully invariant $%
\mathcal{C}$-subcomodule. By Proposition \ref{ro-inn}, a sufficient
condition for $L$ to be fully coprime \emph{(}fully cosemiprime\emph{)} is
that $L$ is fully $M$-coprime \emph{(}fully $M$-cosemiprime\emph{)}, where
the later is also necessary in case $M$ is self-injective.
\end{remark}

\begin{lemma}
\label{inn-ideal}Let $X,Y\subseteq M$ be any $R$-submodules. Then%
\begin{equation}
(X:_{M}^{\mathcal{C}}Y)\subseteq \mathrm{Ke}(\mathrm{An}(X)\circ ^{op}%
\mathrm{An}(Y)),
\end{equation}%
with equality in case $M$ is self-cogenerator and $Y\subseteq M$ is a $%
\mathcal{C}$-subcomodule.
\end{lemma}

\begin{Beweis}
Let $m\in (X:_{M}^{\mathcal{C}}Y)$ be arbitrary. Then for all $f\in \mathrm{%
An}(X)$ we have $f(m)=y$ for some $y\in Y$ and so for each $g\in \mathrm{An}%
(Y)$ we get
\begin{equation*}
(f\circ ^{op}g)(m)=(g\circ f)(m)=g(f(m))=g(y)=0,
\end{equation*}%
i.e. $(X:_{M}^{\mathcal{C}}Y)\subseteq \mathrm{Ke}(\mathrm{An}(X)\circ ^{op}%
\mathrm{An}(Y)).$

Assume now that $M$ is self-cogenerator and that $Y\subseteq M$ is a $%
\mathcal{C}$-subcomodule (so that $\mathrm{KeAn}(Y)=Y$ by \ref{An-Ke} (2)).
If $m\in \mathrm{Ke}(\mathrm{An}(X)\circ ^{op}\mathrm{An}(Y))$ and $f\in
\mathrm{An}(X)$ are arbitrary, then by our choice%
\begin{equation*}
g(f(m))=(f\circ ^{op}g)(m)=0\text{ for all }g\in \mathrm{An}(Y),
\end{equation*}%
so $f(m)\in \mathrm{KeAn}(Y)=Y,$ i.e. $m\in (X:_{M}^{\mathcal{C}}Y).$ Hence,
$(X:_{M}^{\mathcal{C}}Y)=\mathrm{Ke}(\mathrm{An}(X)\circ ^{op}\mathrm{An}%
(Y)).\blacksquare $
\end{Beweis}

\begin{proposition}
\label{corad=}Let $M$ be self-cogenerator. Then
\begin{equation*}
\mathrm{EP}(M)\subseteq \mathrm{CPSpec}(M)\text{ and }\mathrm{ESP}%
(M)\subseteq \mathrm{CSP}(M)
\end{equation*}%
with equality, if $M$ is intrinsically injective self-cogenerator \emph{(}%
whence $\mathrm{EPcorad}(M)=\mathrm{CPcorad}(M)$\emph{)}.
\end{proposition}

\begin{Beweis}
Assume $M$ to be self-cogenerator. Let $0\neq K\subseteq M$ be a fully
invariant $\mathcal{C}$-subcomodule that is $\mathrm{E}$-prime ($\mathrm{E}$%
-semiprime) in $M,$ and suppose $X,$ $Y\subseteq M$ are two (equal) fully
invariant $\mathcal{C}$-subcomodules with $K\subseteq (X:_{M}^{\mathcal{C}%
}Y).$ Then we have by Lemma \ref{inn-ideal} (1)%
\begin{equation*}
\mathrm{An}(X)\circ ^{op}\mathrm{An}(Y)\subseteq \mathrm{AnKe}(\mathrm{An}%
(X)\circ ^{op}\mathrm{An}(Y))\subseteq \mathrm{An}(X:_{M}^{\mathcal{C}%
}Y)\subseteq \mathrm{An}(K).
\end{equation*}%
By assumption $\mathrm{An}(K)\vartriangleleft \mathrm{E}_{M}^{\mathcal{C}}$
is prime (semiprime), hence $\mathrm{An}(X)\subseteq \mathrm{An}(K),$ so
that $K=\mathrm{KeAn}(K)\subseteq \mathrm{KeAn}(X)=X;$ or $\mathrm{An}%
(Y)\subseteq \mathrm{An}(K)$ so that $K=\mathrm{KeAn}(K)\subseteq \mathrm{%
KeAn}(Y)=Y.$ Consequently $K$ is fully $M$-coprime (fully $M$-cosemiprime).

Assume now that $M$ is intrinsically injective self-cogenerator. Let $0\neq
K\subseteq M$ be an fully $M$-coprime (fully $M$-cosemiprime) $\mathcal{C}$%
-subcomodule and consider the proper two-sided ideal $\mathrm{An}%
(K)\vartriangleleft \mathrm{E}_{M}^{\mathcal{C}}.$ Suppose $%
I,J\vartriangleleft \mathrm{E}_{M}^{\mathcal{C}}$ are two (equal) ideals
with $I\circ ^{op}J\subseteq \mathrm{An}(K)$ and $I_{\mathrm{E}_{M}^{%
\mathcal{C}}},$ $J_{\mathrm{E}_{M}^{\mathcal{C}}}$ are \emph{finitely
generated. }Consider the fully invariant $\mathcal{C}$-subcomodules $X:=%
\mathrm{Ke}(I),$ $Y:=\mathrm{Ke}(J)$ of $M.$ Since $M$ is self-cogenerator,
it follows by Lemma \ref{inn-ideal} that%
\begin{equation*}
K=\mathrm{KeAn}(K)\subseteq \mathrm{Ke}(I\circ ^{op}J)=\mathrm{Ke}(\mathrm{An%
}(X)\circ ^{op}\mathrm{An}(Y))=(X:_{M}^{\mathcal{C}}Y).
\end{equation*}%
Since $K$ is fully $M$-coprime (fully $M$-cosemiprime), we conclude that $%
K\subseteq X$ so that $I=\mathrm{AnKe}(I)=\mathrm{An}(X)\subseteq \mathrm{An}%
(K);$ or $K\subseteq Y$ so that $J=\mathrm{AnKe}(J)=\mathrm{An}(Y)\subseteq
\mathrm{An}(K).$ Consequently $\mathrm{An}(K)\vartriangleleft \mathrm{E}%
_{M}^{\mathcal{C}}$ is prime (semiprime), i.e. $K$ is $\mathrm{E}$-prime ($%
\mathrm{E}$-semiprime) in $M.\blacksquare $
\end{Beweis}

\begin{remark}
\label{coprime-tau}It follows from Proposition \ref{corad=} that a
sufficient condition for $M$ to be fully coprime (fully cosemiprime),\emph{\
}in case $M$ is self-cogenerator, is that $M$ is $\mathrm{E}$-prime ($%
\mathrm{E}$-semiprime), where the later is also necessary in case $M$ is
intrinsically injective self-cogenerator.
\end{remark}

\qquad As a direct consequence of Propositions \ref{P-rad}, \ref{corad=} we
have

\begin{proposition}
\label{Prad=CPcorad}Let $M$ be intrinsically injective self-cogenerator and $%
\mathrm{E}_{M}^{\mathcal{C}}$ be right Noetherian. Then%
\begin{equation}
\mathrm{Prad}(\mathrm{E}_{M}^{\mathcal{C}})=\mathrm{An}(\mathrm{CPcorad}(M))%
\text{ and }\mathrm{CPcorad}(M)=\mathrm{Ke}(\mathrm{Prad}(\mathrm{E}_{M}^{%
\mathcal{C}})).  \label{CPcorad=}
\end{equation}
\end{proposition}

\qquad Using Proposition \ref{Prad=CPcorad}, a similar proof to that of
Corollary \ref{ME-semi} yields:

\begin{corollary}
Let $M$ be intrinsically injective self-cogenerator and $\mathrm{E}_{M}^{%
\mathcal{C}}$ be right Noetherian. Then%
\begin{equation*}
M\text{ is fully cosemiprime}\Leftrightarrow M=\mathrm{CPcorad}(M).
\end{equation*}
\end{corollary}

\begin{corollary}
Let $_{A}\mathcal{C}$ be locally projective and $M$ be self injective
self-cogenerator. If $M$ is Artinian \emph{(}e.g. $A$ is right Artinian and $%
M$ is finitely generated\emph{)}, then

\begin{enumerate}
\item $\mathrm{Prad}(\mathrm{E}_{M}^{\mathcal{C}})=\mathrm{An}(\mathrm{%
CPcorad}(M))$ and $\mathrm{CPcorad}(M)=\mathrm{Ke}(\mathrm{Prad}(\mathrm{E}%
_{M}^{\mathcal{C}}).$

\item $M$ is fully cosemiprime $\Leftrightarrow M=\mathrm{CPcorad}(M).$
\end{enumerate}
\end{corollary}

\subsection*{Comodules with rings of colinear endomorphisms right Artinian}

\qquad Under the \emph{strong assumption} $\mathrm{E}_{M}^{\mathcal{C}}:=%
\mathrm{End}^{\mathcal{C}}(M)^{op}$ is right Artinian, several primeness and
coprimeness properties of the non-zero right $\mathcal{C}$-comodule $M$
coincide and become, in case $_{A}\mathcal{C}$ locally projective,
equivalent to $M$ being simple as a $(^{\ast }\mathcal{C},\mathrm{E}_{M}^{%
\mathcal{C}})$-bimodule. This follows from the fact that right Artinian
prime rings are simple.

\qquad If $M$ has no non-trivial fully invariant $\mathcal{C}$-subcomodules,
then it is obviously fully coprime. The following result gives a partial
converse:

\begin{theorem}
\label{E-Artinian}Let $M$ be intrinsically injective self-cogenerator and
assume $\mathrm{E}_{M}^{\mathcal{C}}$ to be right Artinian. Then the
following are equivalent:

\begin{enumerate}
\item $M$ is $\mathrm{E}$-prime \emph{(}i.e. $\mathrm{E}_{M}^{\mathcal{C}}$
is a prime ring\emph{)};

\item $\mathrm{E}_{M}^{\mathcal{C}}$ is simple;

\item $M$ has no non-trivial fully invariant $\mathcal{C}$-subcomodules;

\item $M$ is fully coprime.
\end{enumerate}
\end{theorem}

\begin{Beweis}
Let $M$ be intrinsically injective self-cogenerator and assume $\mathrm{E}%
_{M}^{\mathcal{C}}$ to be right Artinian.

$(1)\Rightarrow (2):$ Right Artinian prime rings are simple (e.g. \cite[4.5
(2)]{Wis91}).

$(2)\Rightarrow (3):$ Since $M$ is self-cogenerator, this follows by Theorem %
\ref{E-simple}.

$(3)\Rightarrow (4):$ Trivial.

$(4)\Rightarrow (1):$ Since $M$ be intrinsically injective self-cogenerator,
this follows by \ref{coprime-tau}.$\blacksquare $
\end{Beweis}

\begin{proposition}
\label{coprime=simple}Let $_{A}\mathcal{C}$ be locally projective and $M$ be
self-injective self-cogenerator. If any of the following additional
conditions is satisfied, then $M$ is fully coprime if and only if $M$ is
simple as a $(^{\ast }\mathcal{C},\mathrm{E}_{M}^{\mathcal{C}})$-bimodule:

\begin{enumerate}
\item $M$ has finite length; \emph{or}

\item $A$ is right Artinian and $M_{A}$ is finitely generated; \emph{or}

\item $M$ is Artinian and self-projective.
\end{enumerate}
\end{proposition}

\begin{Beweis}
By Theorem \ref{E-Artinian}, it suffices to show that $\mathrm{E}_{M}^{%
\mathcal{C}}=\mathrm{End}(_{^{\ast }\mathcal{C}}M)^{op}$ is right Artinian
under each of the additional conditions.

\begin{enumerate}
\item By assumption $M$ is self-injective and Artinian (semi-injective and
Noetherian) and it follows then by Proposition \ref{End} that $\mathrm{E}%
_{M}^{\mathcal{C}}$ is right Noetherian (semiprimary). Applying Hopkins
Theorem (e.g. \cite[31.4]{Wis91}), we conclude that $\mathrm{E}_{M}^{%
\mathcal{C}}$ is right Artinian.

\item If $A$ is right Artinian and $_{A}\mathcal{C}$ is locally projective,
then every finitely generated right $\mathcal{C}$-comodule has finite length
by \cite[Corollary 2.25]{Abu03}.

\item Since $M$ is Artinian, self-injective and self-projective, $\mathrm{E}%
_{M}^{\mathcal{C}}$ is right Artinian by Proposition \ref{End} (2).$%
\blacksquare $
\end{enumerate}
\end{Beweis}

\subsection*{Fully coprimeness versus irreducibility}

\qquad In what follows we clarify, under suitable conditions, the relation
between fully coprime and irreducible comodules:

\begin{proposition}
\label{ro-union}Let $\{K_{\lambda }\}_{\Lambda }$ be a family of non-zero
fully invariant $\mathcal{C}$-subcomodules of $M,$ such that for any $\gamma
,\delta \in \Lambda $ either $K_{\gamma }$ $\subseteq K_{\delta }$ or $%
K_{\delta }\subseteq K_{\gamma },$ and consider the fully invariant $%
\mathcal{C}$-subcomodule $K:=\sum\limits_{\lambda \in \Lambda }K_{\lambda
}=\bigcup\limits_{\lambda \in \Lambda }K_{\lambda }\subseteq M.$ If $%
K_{\lambda }\in \mathrm{CPSpec}(M)$ for all $\lambda \in \Lambda ,$ then $%
K\in \mathrm{CPSpec}(M).$

\begin{Beweis}
Let $X,Y\subseteq M$ be any fully invariant $\mathcal{C}$-subcomodules with $%
K\subseteq (X:_{M}^{\mathcal{C}}Y)$ and suppose $K\nsubseteq X.$ We \textbf{%
claim} that $K\subseteq Y.$

Since $K\nsubseteq X,$ there exists some $\lambda _{0}\in \Lambda $ with $%
K_{\lambda _{0}}\nsubseteq X.$ Since $K_{\lambda _{0}}\subseteq (X:_{M}^{%
\mathcal{C}}Y),$ it follows from the assumption $K_{\lambda _{0}}\in \mathrm{%
CPSpec}(M)$ that $K_{\lambda _{0}}\subseteq Y.$ Let $\lambda \in \Lambda $
be arbitrary. If $K_{\lambda }\subseteq K_{\lambda _{0}},$ then $K_{\lambda
}\subseteq Y.$ If otherwise $K_{\lambda _{0}}\subseteq K_{\lambda },$ then
the inclusion $K_{\lambda }\subseteq (X:_{M}^{\mathcal{C}}Y)$ implies $%
K_{\lambda }\subseteq Y$ (since $K_{\lambda }\subseteq X$ would imply $%
K_{\lambda _{0}}\subseteq X,$ a contradiction). So $K:=\bigcup\limits_{%
\lambda \in \Lambda }K_{\lambda }\subseteq Y.\blacksquare $
\end{Beweis}
\end{proposition}

\begin{corollary}
\label{M-ro-union}Let $M=\sum\limits_{\lambda \in \Lambda }M_{\lambda },$
where $\{M_{\lambda }\}_{\Lambda }$ is a family of non-zero fully invariant $%
\mathcal{C}$-subcomodules of $M$ such that for any $\gamma ,\delta \in
\Lambda $ either $M_{\gamma }$ $\subseteq M_{\delta }$ or $M_{\delta
}\subseteq M_{\gamma }.$ If $M_{\lambda }\in \mathrm{CPSpec}(M)$ for each $%
\lambda \in \Lambda ,$ then $M$ is fully coprime.
\end{corollary}

\begin{proposition}
\label{ro-irr}Let $0\neq K\subseteq M$ be a non-zero fully invariant $%
\mathcal{C}$-subcomodule. If $K\in \mathrm{CPSpec}(M),$ then $K$ has no
decomposition as an internal direct sum of non-trivial fully invariant $%
\mathcal{C}$-subcomodules.
\end{proposition}

\begin{Beweis}
Let $K\in \mathrm{CPSpec}(M)$ and suppose $K:=K_{\lambda _{0}}\oplus
\sum\limits_{\lambda \neq \lambda _{0}}K_{\lambda },$ an internal direct sum
of non-trivial fully invariant $\mathcal{C}$-subcomodules. Then $K\subseteq
(K_{\lambda _{0}}:_{M}^{\mathcal{C}}\sum\limits_{\lambda \neq \lambda
_{0}}K_{\lambda })$ and it follows that $K\subseteq K_{\lambda _{0}}$ or $%
K\subseteq \sum\limits_{\lambda \neq \lambda _{0}}K_{\lambda }$
(contradiction).$\blacksquare $
\end{Beweis}

\begin{corollary}
\label{M-pre-ind}If $M$ is fully coprime, then $M$ has no decomposition as
an internal direct sum of non-trivial fully invariant $\mathcal{C}$%
-subcomodules.
\end{corollary}

As a direct consequence of Corollary \ref{M-pre-ind} we get a restatement of
Theorem \ref{E-coprime-irr}:

\begin{theorem}
\label{coprime-irr}Let $_{A}\mathcal{C}$ locally projective and $M$ be
self-injective self-cogenerator with $\mathrm{End}^{\mathcal{C}}(M)$
commutative. If $M$ is fully coprime, then $M$ is irreducible.
\end{theorem}

\section{Primeness and Coprimeness Conditions for Corings}

\qquad Throughout this section $(\mathcal{C},\Delta ,\varepsilon )$ is a
non-zero coring. We consider in what follows several coprimeness
(cosemiprimeness) and primeness (semiprimeness) properties of $\mathcal{C},$
considered as an object in the category $\mathbb{M}^{\mathcal{C}}$ of right $%
\mathcal{C}$-comodules, denoted by $\mathcal{C}^{r},$ as well as an object
in the category $^{\mathcal{C}}\mathbb{M}$ of left $\mathcal{C}$-comodules,
denoted by $\mathcal{C}^{l}.$ In particular, we clarify the relation between
these properties and the simplicity (semisimplicity) of $\mathcal{C}.$
Several results in this section can be obtained directly from the
corresponding ones in the previous sections. Moreover, we state many of
these in the case $A$ is a QF ring, as in this case $\mathcal{C}$ is an
injective cogenerator in both the categories of right and left $\mathcal{C}$%
-comodules by Lemma \ref{MotC-inj}.

\begin{punto}
\label{An=per}(e.g. \cite[17.8.]{BW03}) We have an isomorphism of $R$%
-algebras%
\begin{equation*}
\phi _{r}:\mathcal{C}^{\ast }\rightarrow \mathrm{End}^{\mathcal{C}}(\mathcal{%
C})^{op},\text{ }f\mapsto \lbrack c\mapsto c\leftharpoonup f:=\sum
f(c_{1})c_{2}]
\end{equation*}%
with inverse map $\psi _{r}:g\mapsto \varepsilon \circ g$ $,$ and there is a
ring morphism $\iota _{r}:A\longrightarrow (\mathcal{C}^{\ast })^{op},$ $%
a\mapsto \varepsilon (a-).$

Similarly, we have an isomorphism of $R$-algebras
\begin{equation*}
\phi _{l}:\text{ }^{\ast }\mathcal{C}\rightarrow \text{ }^{\mathcal{C}}%
\mathrm{End}(\mathcal{C}),\text{ }f\mapsto \lbrack c\mapsto f\rightharpoonup
c:=\sum c_{1}f(c_{2})]
\end{equation*}%
with inverse map $\psi _{l}:g\mapsto \varepsilon \circ g,$ and there is a
ring morphism $\iota _{l}:A\longrightarrow (^{\ast }\mathcal{C})^{op},$ $%
a\mapsto \varepsilon (-a).$
\end{punto}

\begin{definition}
\begin{enumerate}
\item We call a right (left) $A$-submodule $K\subseteq \mathcal{C}$ a \emph{%
right }(\emph{left})\emph{\ }$\mathcal{C}$\emph{-coideal}, iff $K$ is a
right (left) $\mathcal{C}$-subcomodule of $\mathcal{C}$ with structure map
the restriction of $\Delta _{\mathcal{C}}$ to $K.$

\item We call an $(A,A)$-subbimodule $B\subseteq \mathcal{C}$ a $\mathcal{C}$%
\emph{-bicoideal, }iff $B$ is a $\mathcal{C}$-subbicomodule of $\mathcal{C}$
with structure map the restriction of $\Delta _{\mathcal{C}}$ to $B;$

\item We call an $(A,A)$-subbimodule $\mathcal{D}\subseteq \mathcal{C}$ an $%
A $\emph{-subcoring, }iff $\mathcal{D}$ is an $A$-coring with structure maps
the restrictions of $\Delta _{\mathcal{C}}$ and $\varepsilon _{\mathcal{C}}$
to $\mathcal{D}.$
\end{enumerate}
\end{definition}

\begin{notation}
With $\mathcal{R}(\mathcal{C})$ ($\mathcal{R}_{f.i.}(\mathcal{C})$) we
denote the class of (fully invariant) right $\mathcal{C}$-coideals and with $%
\mathcal{I}_{r}(\mathcal{C}^{\ast })$ ($\mathcal{I}_{t.s.}(\mathcal{C}^{\ast
})$) the class of right (two-sided) ideals of $\mathcal{C}^{\ast }.$
Analogously, we denote with $\mathcal{L}(\mathcal{C})$ ($\mathcal{L}_{f.i.}(%
\mathcal{C})$) the class of (fully invariant) left $\mathcal{C}$-coideals
and with $\mathcal{I}_{l}(^{\ast }\mathcal{C})$ ($\mathcal{I}_{t.s.}(^{\ast }%
\mathcal{C})$) the class of left (two-sided) ideals of $^{\ast }\mathcal{C}.$
With $\mathcal{B}(\mathcal{C})$ we denote the class of $\mathcal{C}$%
-bicoideals and for each $B\in \mathcal{B}(\mathcal{C})$ we write $B^{r}$ ($%
B^{l}$) to indicate that we consider $B$ as an object in the category of
right (left) $\mathcal{C}$-comodules.
\end{notation}

\begin{remarks}
\label{coideal-rem}For $\varnothing \neq I\subseteq \mathcal{C}^{\ast }$ ($%
\varnothing \neq I\subseteq $ $^{\ast }\mathcal{C}$) and $\varnothing \neq
K\subseteq \mathcal{C},$ set%
\begin{equation*}
I^{\bot (\mathcal{C})}:=\bigcap\limits_{f\in I}\{c\in \mathcal{C}\mid
f(c)=0\}
\end{equation*}%
and%
\begin{equation*}
\begin{tabular}{lll}
$K^{\bot (^{\ast }\mathcal{C})}:=\{f\in $ $^{\ast }\mathcal{C}\mid f(K)=0\};$
&  & $K^{\bot (\mathcal{C}^{\ast })}:=\{f\in \mathcal{C}^{\ast }\mid
f(K)=0\}.$%
\end{tabular}%
\end{equation*}

\begin{enumerate}
\item If $_{A}\mathcal{C}$ is flat, then a right $A$-submodule $K\subseteq
\mathcal{C}$ is a right $\mathcal{C}$-coideal, iff $\Delta (K)\subseteq
K\otimes _{A}\mathcal{C}.$

If $\mathcal{C}_{A}$ is flat, then a left $A$-submodule $K\subseteq \mathcal{%
C}$ is a left $\mathcal{C}$-coideal, iff $\Delta (K)\subseteq \mathcal{C}%
\otimes _{A}K.$

If $_{A}\mathcal{C}$ and $\mathcal{C}_{A}$ are flat, then an $A$-subbimodule
$B\subseteq \mathcal{C}$ is a $\mathcal{C}$-bicoideal, iff $\Delta
(B)\subseteq (B\otimes _{A}\mathcal{C})\cap (\mathcal{C}\otimes _{A}B).$

If $_{A}\mathcal{C}$ and $\mathcal{C}_{A}$ are flat, then an $A$-subbimodule
$\mathcal{D}\subseteq \mathcal{C}$ is a subcoring, iff $\Delta (\mathcal{D}%
)\subseteq \mathcal{D}\otimes _{A}\mathcal{D}.$

\item Every $A$-subcoring $\mathcal{D}\subseteq \mathcal{C}$ is a $\mathcal{C%
}$-bicoideal in the canonical way.

If $B\subseteq \mathcal{C}$ is a $\mathcal{C}$-bicoideal that is pure as a
left and as a right $A$-submodule, then we have by \cite[40.16]{BW03}:%
\begin{equation*}
\Delta (B)\subseteq (B\otimes _{A}\mathcal{C})\cap (\mathcal{C}\otimes
_{A}B)=B\otimes _{A}B,
\end{equation*}%
i.e. $B\subseteq \mathcal{C}$ is an $A$-subcoring.

\item If $\mathcal{C}_{A}$ (respectively $_{A}\mathcal{C}$) is locally
projective, then $\mathcal{R}_{f.i.}(\mathcal{C})=\mathcal{B}(\mathcal{C})$
(respectively $\mathcal{L}_{f.i.}(\mathcal{C})=\mathcal{B}(\mathcal{C})$):
if $B\subseteq \mathcal{C}$ is a fully invariant right (left) $\mathcal{C}$%
-coideal, then $B\subseteq \mathcal{C}$ is a right $\mathcal{C}^{\ast }$%
-submodule (left $^{\ast }\mathcal{C}$-submodule) and it follows by
Proposition \ref{comod=sg} that $B\subseteq \mathcal{C}$ is a $\mathcal{C}$%
-subbicomodule with structure map the restriction of $\Delta _{\mathcal{C}}$
to $B,$ i.e. $B$ is a $\mathcal{C}$-bicoideal.

\item Let $\mathcal{C}_{A}$ ($_{A}\mathcal{C}$) be locally projective. If $%
P\vartriangleleft \mathcal{C}^{\ast }$ ($P\vartriangleleft $ $^{\ast }%
\mathcal{C}$) is a two-sided ideal, then the fully invariant right (left) $%
\mathcal{C}$-coideal $B:=\mathrm{ann}_{\mathcal{C}}(P)\subseteq \mathcal{C}$
is a $\mathcal{C}$-bicoideal.

\item If $_{A}\mathcal{C}$ is locally projective and $I\vartriangleleft _{r}$
$^{\ast }\mathcal{C}$ is a right ideal, then the left $^{\ast }\mathcal{C}$%
-submodule $I^{\bot (\mathcal{C})}\subseteq \mathcal{C}$ is a right $%
\mathcal{C}$-coideal.

If $\mathcal{C}_{A}$ is locally projective and $I\vartriangleleft _{l}%
\mathcal{C}^{\ast }$ is a left ideal, then the right $\mathcal{C}^{\ast }$%
-submodule $I^{\bot (\mathcal{C})}\subseteq \mathcal{C}$ is a left $\mathcal{%
C}$-coideal.

\item If $K\subseteq \mathcal{C}$ is a (fully invariant) right $\mathcal{C}$%
-coideal, then $K^{\bot (\mathcal{C}^{\ast })}=\mathrm{ann}_{\mathcal{C}%
^{\ast }}(K)\simeq \mathrm{An}_{\mathrm{E}_{\mathcal{C}}^{\mathcal{C}}}(K);$
in particular $K^{\bot (\mathcal{C}^{\ast })}\subseteq \mathcal{C}^{\ast }$
is a right (two-sided) ideal.

If $K\subseteq \mathcal{C}$ is a (fully invariant) left $\mathcal{C}$%
-coideal, then $K^{\bot (^{\ast }\mathcal{C})}=\mathrm{ann}_{^{\ast }%
\mathcal{C}}(K)\simeq \mathrm{An}_{_{\mathcal{C}}^{\mathcal{C}}\mathrm{E}%
}(K);$ in particular $K^{\bot (^{\ast }\mathcal{C})}\subseteq $ $^{\ast }%
\mathcal{C}$ is a left (two-sided) ideal.

\item If $A_{A}$ is an injective cogenerator and $_{A}\mathcal{C}$ is flat,
then for every right ideal $I\vartriangleleft _{r}\mathcal{C}^{\ast }$ we
have $\mathrm{ann}_{\mathcal{C}}(I)=I^{\bot (\mathcal{C})}:$ Write $%
I=\bigcup\limits_{\lambda \in \Lambda }I_{\lambda },$ where $I_{\lambda
}\vartriangleleft _{r}\mathcal{C}^{\ast }$ is a finitely generated right
ideal for each $\lambda \in \Lambda .$ If $\mathrm{ann}_{\mathcal{C}%
}(I_{\lambda _{0}})\varsubsetneqq I_{\lambda _{0}}^{\bot (\mathcal{C})}$ for
some $\lambda _{0}\in \Lambda ,$ then $\mathrm{Hom}_{A}(\mathcal{C}/\mathrm{%
ann}_{\mathcal{C}}(I_{\lambda _{0}}),A)\nsubseteq \mathrm{Hom}_{A}(\mathcal{C%
}/I_{\lambda _{0}}^{\bot (\mathcal{C})},A)$ (since $A_{A}$ is a
cogenerator). Since $A_{A}$ is injective, $\mathcal{C}$ is injective in $%
\mathbb{M}^{\mathcal{C}}$ by Lemma \ref{MotC-inj} and it follows by \ref%
{An-Ke} (3-b) and the remarks above that $I_{\lambda _{0}}=\mathrm{ann}_{%
\mathcal{C}^{\ast }}(\mathrm{ann}_{\mathcal{C}}(I_{\lambda _{0}}))=(\mathrm{%
ann}_{\mathcal{C}}(I_{\lambda _{0}}))^{\bot (\mathcal{C}^{\ast
})}\nsubseteqq I_{\lambda _{0}}^{\bot (\mathcal{C})\bot (\mathcal{C}^{\ast
})}$ (a contradiction). So $\mathrm{ann}_{\mathcal{C}}(I_{\lambda
})=I_{\lambda }^{\bot (\mathcal{C})}$ for each $\lambda \in \Lambda $ and we
get%
\begin{equation*}
\mathrm{ann}_{\mathcal{C}}(I)=\bigcap_{\lambda \in \Lambda }\mathrm{ann}_{%
\mathcal{C}}(I_{\lambda })=\bigcap_{\lambda \in \Lambda }I_{\lambda }^{\bot (%
\mathcal{C})}=(\bigcup_{\lambda \in \Lambda }I_{\lambda })^{\bot (\mathcal{C}%
)}=I^{\bot (\mathcal{C})}.\blacksquare
\end{equation*}
\end{enumerate}
\end{remarks}

\subsection*{Sufficient and necessary conditions}

\qquad The following result gives sufficient and necessary conditions for
the dual rings of the non-zero coring $\mathcal{C}$ to be prime
(respectively semiprime, domain, reduced) generalizing results of \cite%
{YDZ90} for coalgebras over base fields. Several results follow directly
from previous sections recalling the isomorphism of rings $\mathcal{C}^{\ast
}\simeq \mathrm{End}^{\mathcal{C}}(\mathcal{C})^{op}.$ Analogous statements
can be formulated for $^{\ast }\mathcal{C}\simeq $ $^{\mathcal{C}}\mathrm{End%
}(\mathcal{C}).$

\qquad Theorem \ref{iff-M} yields directly.

\begin{theorem}
\label{sn-cor}Let $_{A}\mathcal{C}$ be flat.

\begin{enumerate}
\item $\mathcal{C}^{\ast }$ is prime \emph{(}domain\emph{),} if $\mathcal{C}$
$=\mathcal{C}f\mathcal{C}^{\ast }$ \emph{(}$\mathcal{C}=\mathcal{C}f$\emph{)}
for all $0\neq f\in \mathcal{C}^{\ast }.$ If $\mathcal{C}$ is coretractable
in $\mathbb{M}^{\mathcal{C}},$ then $\mathcal{C}^{\ast }$ is prime \emph{(}%
domain\emph{)} if and only if $\mathcal{C}=\mathcal{C}f\mathcal{C}^{\ast }$
\emph{(}$\mathcal{C}=\mathcal{C}f$\emph{)} $\forall $ $0\neq f\in \mathcal{C}%
^{\ast }.$

\item $\mathcal{C}^{\ast }$ is semiprime \emph{(}reduced\emph{),} if $%
\mathcal{C}f=\mathcal{C}f\mathcal{C}^{\ast }f$ \emph{(}$\mathcal{C}f=%
\mathcal{C}f^{2}$\emph{)} for all $0\neq f\in \mathcal{C}^{\ast }.$ If $%
\mathcal{C}$ is self-cogenerator in $\mathbb{M}^{\mathcal{C}}\emph{,}$ then $%
\mathcal{C}^{\ast }$ is semiprime \emph{(}reduced\emph{)} if and only if $%
\mathcal{C}f=\mathcal{C}f\mathcal{C}^{\ast }f$ \emph{(}$\mathcal{C}f=%
\mathcal{C}f^{2}$\emph{)} $\forall $ $0\neq f\in \mathcal{C}^{\ast }.$
\end{enumerate}
\end{theorem}

\begin{proposition}
\label{r-l-prime}Let $_{A}\mathcal{C}$ and $\mathcal{C}_{A}$ be flat.

\begin{enumerate}
\item Let $\mathcal{C}$ be coretractable in $\mathbb{M}^{\mathcal{C}}$ and $%
\mathcal{C}_{\mathcal{C}^{\ast }}$ satisfy condition \emph{(}**\emph{)}. If $%
\mathcal{C}^{\ast }$ is prime \emph{(}domain\emph{)},\emph{\ }then $^{\ast }%
\mathcal{C}$ is prime \emph{(}domain\emph{)}.

\item Let $\mathcal{C}$ be coretractable in $\mathbb{M}^{\mathcal{C}},$ $^{%
\mathcal{C}}\mathbb{M}$ and $\mathcal{C}_{\mathcal{C}^{\ast }},_{^{\ast }%
\mathcal{C}}\mathcal{C}$ satisfy condition \emph{(}**\emph{)}. Then $%
\mathcal{C}^{\ast }$ is prime \emph{(}domain\emph{)} if and only if\emph{\ }$%
^{\ast }\mathcal{C}$ is prime \emph{(}domain\emph{)}.

\item Let $\mathcal{C}$ be coretractable in $\mathbb{M}^{\mathcal{C}}$ and $%
_{A}\mathcal{C}$ be locally projective. If $\mathcal{C}^{\ast }$ is prime,%
\emph{\ }then $^{\ast }\mathcal{C}$ is prime.

\item Let $\mathcal{C}$ be coretractable in $\mathbb{M}^{\mathcal{C}},$ $^{%
\mathcal{C}}\mathbb{M}$ and $_{A}\mathcal{C},$ $\mathcal{C}_{A}$ be locally
projective. Then $\mathcal{C}^{\ast }$ is prime if and only if $^{\ast }%
\mathcal{C}$ is prime.
\end{enumerate}
\end{proposition}

\begin{Beweis}
\begin{enumerate}
\item Let $\mathcal{C}^{\ast }$ be prime (domain). If $^{\ast }\mathcal{C}$
were not prime (not a domain), then there exists by an analogous statement
of Theorem \ref{sn-cor} some $0\neq f\in $ $^{\ast }\mathcal{C}$ with $%
^{\ast }\mathcal{C}f\mathcal{C}\varsubsetneqq \mathcal{C\ }$($f\mathcal{C}%
\varsubsetneqq \mathcal{C}$). By assumption $\mathcal{C}_{\mathcal{C}^{\ast
}}$ satisfies condition (**) and so there exists some $0\neq h\in \mathcal{C}%
^{\ast }$ such that $(^{\ast }\mathcal{C}f\mathcal{C})h=0$ ($(f\mathcal{C}%
)h=0$). But this implies $\mathcal{C}\neq \mathcal{C}h\mathcal{C}^{\ast }$ ($%
\mathcal{C}\neq \mathcal{C}h$): otherwise $f\mathcal{C}=f(\mathcal{C}h%
\mathcal{C}^{\ast })=((f\mathcal{C})h)\mathcal{C}^{\ast }=0$ ($f\mathcal{C}%
=f(\mathcal{C}h)=(f\mathcal{C})h=0$), which implies $f=0,$ a contradiction.
Since $\mathcal{C}$ is coretractable in $\mathbb{M}^{\mathcal{C}},$ Theorem %
\ref{sn-cor} (1) implies that $\mathcal{C}^{\ast }$ is not prime (not a
domain), which contradicts our assumptions.

\item Follows from (1) by symmetry.

\item The proof is similar to that of (1) recalling that, in case $_{A}%
\mathcal{C}$ locally projective, for any $f\in $ $^{\ast }\mathcal{C},$ the
left $^{\ast }\mathcal{C}$-submodule $^{\ast }\mathcal{C}f\mathcal{C}%
\subseteq \mathcal{C}$ is a right $\mathcal{C}$-subcomodule.

\item Follows from (3) by symmetry.$\blacksquare $
\end{enumerate}
\end{Beweis}

\subsection*{\textrm{E}-Prime versus simple}

\qquad In what follows we show that \textrm{E}-prime corings generalize
simple corings. The results are obtained by direct application of the
corresponding results in the Section 3.

\qquad As a direct consequence of Theorems \ref{r-simple-E} and \ref%
{E-simple} we get

\begin{theorem}
\label{C-r-simple-E}Let $A$ be a QF ring and assume $_{A}\mathcal{C}$ to be
\emph{(}locally\emph{) }projective.

\begin{enumerate}
\item $\mathcal{C}^{r}$ is simple if and only if $\mathcal{C}^{\ast }$ is
right simple.

\item If $\mathcal{C}^{\ast }$ is simple, then $\mathcal{C}$ is simple \emph{%
(}as a $(^{\ast }\mathcal{C},\mathcal{C}^{\ast })$-bimodule\emph{)}.

\item Let $\mathcal{C}^{\ast }$ be right Noetherian. Then $\mathcal{C}^{\ast
}$ is simple if and only if $\mathcal{C}$ is simple \emph{(}as a $(^{\ast }%
\mathcal{C},\mathcal{C}^{\ast })$-bimodule\emph{)}.
\end{enumerate}
\end{theorem}

\begin{corollary}
Let $A$ be a QF ring, $_{A}\mathcal{C},$ $\mathcal{C}_{A}$ be locally
projective, $^{\ast }\mathcal{C}$ be left Noetherian and $\mathcal{C}^{\ast
} $ be right Noetherian. Then%
\begin{equation*}
\mathcal{C}^{\ast }\text{ is simple}\Leftrightarrow \mathcal{C}\text{ is
simple \emph{(}as a }(^{\ast }\mathcal{C},\mathcal{C}^{\ast })\text{-bimodule%
\emph{)} }\Leftrightarrow \text{ }^{\ast }\mathcal{C}\text{ is simple.}
\end{equation*}
\end{corollary}

\begin{proposition}
\label{C-Jac}Let $A$ be a QF ring. If $_{A}\mathcal{C}$ is \emph{(}locally%
\emph{)} projective, then we have%
\begin{equation*}
\mathrm{Jac}(\mathcal{C}^{\ast })=\mathrm{ann}_{\mathcal{C}^{\ast }}(\mathrm{%
Soc}(\mathcal{C}^{r}))=\mathrm{Soc}(\mathcal{C}^{r})^{\bot (\mathcal{C}%
^{\ast })}\text{ and }\mathrm{Soc}(\mathcal{C}^{r})=\mathrm{Jac}(\mathcal{C}%
^{\ast })^{\bot (\mathcal{C})}.
\end{equation*}
\end{proposition}

\begin{Beweis}
The result in (1) follows from Proposition \ref{Jac} (3) recalling the
isomorphisms of $R$-algebras $\mathcal{C}^{\ast }\simeq \mathrm{End}^{%
\mathcal{C}}(\mathcal{C})^{op}$ and Remarks \ref{coideal-rem} (6) $\&\ $(7).$%
\blacksquare $
\end{Beweis}

\begin{corollary}
\label{C-semisimple-semiprimitive}Let $A$ be a QF ring.

\begin{enumerate}
\item If $_{A}\mathcal{C}$ is \emph{(}locally\emph{)} projective, then $%
\mathcal{C}$ is right semisimple if and only if $\mathcal{C}^{\ast }$ is
semiprimitive.

\item If $_{A}\mathcal{C}$ and $\mathcal{C}_{A}$ are \emph{(}locally\emph{)}
projective, then $\mathcal{C}^{\ast }$ is semiprimitive if and only if $%
^{\ast }\mathcal{C}$ is semiprimitive.
\end{enumerate}
\end{corollary}

\subsection*{The wedge product}

\qquad The \emph{wedge product} of subspaces of a given coalgebra $C$ over a
base field was already defined and investigated in \cite[Section 9]{Swe69}.
In \cite{NT01}, the wedge product of subcoalgebras was used to define \emph{%
fully coprime coalgebras.}

\begin{definition}
We define the \emph{wedge product }of a right $A$-submodule $K\subseteq
\mathcal{C}$ and a left $A$-submodule $L\subseteq \mathcal{C}$ as%
\begin{equation*}
K\wedge L:=\Delta ^{-1}(\func{Im}(K\otimes _{A}\mathcal{C})+\func{Im}(%
\mathcal{C}\otimes _{A}L))=\mathrm{Ker}((\pi _{K}\otimes \pi _{L})\circ
\Delta :\mathcal{C}\longrightarrow \mathcal{C}/K\otimes _{A}\mathcal{C}/L).
\end{equation*}
\end{definition}

\begin{remark}
(\cite[Proposition 9.0.0.]{Swe69}) Let $C$ be a coalgebra over a base field
and $K,L\subseteq C$ be any subspaces. Then $K\wedge L=(K^{\bot (C^{\ast
})}\ast L^{\bot (C^{\ast })})^{\bot (C)}.$ If moreover $K$ is a left $C$%
-coideal and $L$ is a right $C$-coideal, then $K\wedge L\subseteq C$ is a
subcoalgebra.
\end{remark}

\begin{lemma}
\label{kap}\emph{(See \cite[Corollary 2.9.]{Abu03})} Let $K,L\subseteq
\mathcal{C}$ be $A$-subbimodules.

\begin{enumerate}
\item Consider the canonical $A$-bilinear map%
\begin{equation*}
\kappa _{l}:K^{\bot (^{\ast }\mathcal{C})}\otimes _{A}L^{\bot (^{\ast }%
\mathcal{C})}\rightarrow \text{ }^{\ast }(\mathcal{C}\otimes _{A}\mathcal{C}%
),\text{ }[f\otimes _{A}g\mapsto (c\otimes _{A}c^{\prime })=g(cf(c^{\prime
}))].
\end{equation*}%
If $A$ is right Noetherian, $\mathcal{C}_{A}$ is flat and $L^{\perp (^{\ast }%
\mathcal{C})\bot }\subseteq \mathcal{C}$ is pure as a right $A$-module, then%
\begin{equation}
(\kappa _{l}(K^{\bot (^{\ast }\mathcal{C})}\otimes _{A}L^{\bot (^{\ast }%
\mathcal{C})}))^{\perp (\mathcal{C}\otimes _{A}\mathcal{C})}=L^{\bot (^{\ast
}\mathcal{C})\bot }\otimes _{A}\mathcal{C}+\mathcal{C}\otimes _{A}K^{\bot
(^{\ast }\mathcal{C})\bot }.  \label{k-l}
\end{equation}

\item Consider the canonical $A$-bilinear map
\begin{equation*}
\kappa _{r}:L^{\perp (\mathcal{C}^{\ast })}\otimes _{A}K^{\perp (\mathcal{C}%
^{\ast })}\rightarrow (\mathcal{C}\otimes _{A}\mathcal{C})^{\ast },\text{ }%
[g\otimes _{A}f\mapsto (c^{\prime }\otimes _{A}c)=g(f(c^{\prime })c)].
\end{equation*}%
If $A$ is left Noetherian, $_{A}\mathcal{C}$ is flat and $L^{\perp (\mathcal{%
C}^{\ast })\bot }\subseteq \mathcal{C}$ is pure as a left $A$-module, then%
\begin{equation}
(\kappa _{r}(L\otimes _{A}K))^{\perp (\mathcal{C}\otimes _{A}\mathcal{C}%
)}=K^{\perp (\mathcal{C}^{\ast })\bot }\otimes _{A}\mathcal{C}+\mathcal{C}%
\otimes _{A}L^{\perp (\mathcal{C}^{\ast })\bot }.  \label{k-r}
\end{equation}%
\qquad
\end{enumerate}
\end{lemma}

\begin{definition}
For $R$-submodules $K,L\subseteq \mathcal{C}$ we set%
\begin{equation*}
\begin{tabular}{lll}
$(K:_{\mathcal{C}^{r}}L)$ & $:=$ & $\bigcap \{f^{-1}(Y)\mid f\in \mathrm{End}%
^{\mathcal{C}}(\mathcal{C})^{op}$ and $f(K)=0\}$ \\
& $=$ & $\bigcap \{c\in \mathcal{C}\mid c\leftharpoonup f\in L\text{ for all
}f\in \mathrm{ann}_{\mathcal{C}^{\ast }}(K)\}.$%
\end{tabular}%
\end{equation*}%
and%
\begin{equation*}
\begin{tabular}{lll}
$(K:_{\mathcal{C}^{l}}L)$ & $:=$ & $\bigcap \{f^{-1}(L)\mid f\in \text{ }^{%
\mathcal{C}}\mathrm{End}(\mathcal{C})\text{ and }f(K)=0\}$ \\
& $=$ & $\bigcap \{c\in \mathcal{C}\mid f\rightharpoonup c\in L$ for all $%
f\in \mathrm{ann}_{^{\ast }\mathcal{C}}(K)\}.$%
\end{tabular}%
\end{equation*}

If $K,L\subseteq \mathcal{C}$ are right (left) $\mathcal{C}$-coideals, then
we call $(K:_{\mathcal{C}^{r}}L)$ ($(K:_{\mathcal{C}^{l}}L)$) the \emph{%
internal coproduct }of $X$ and $Y$ in $\mathbb{M}^{\mathcal{C}}$ (in $^{%
\mathcal{C}}\mathbb{M}$).
\end{definition}

\begin{lemma}
\label{biprod=}Let $K,L\subseteq \mathcal{C}$ be $\mathcal{C}$-bicoideals.

\begin{enumerate}
\item If $_{A}\mathcal{C}$ is flat and $\mathcal{C}$ is self-cogenerator in $%
\mathbb{M}^{\mathcal{C}},$ then%
\begin{equation*}
(K:_{\mathcal{C}^{r}}L)=\mathrm{ann}_{\mathcal{C}}(\mathrm{ann}_{\mathcal{C}%
^{\ast }}(K)\ast ^{r}\mathrm{ann}_{\mathcal{C}^{\ast }}(L)).
\end{equation*}

\item If $\mathcal{C}_{A}$ is flat and $\mathcal{C}$ is self-cogenerator in $%
^{\mathcal{C}}\mathbb{M},$ then%
\begin{equation*}
(K:_{\mathcal{C}^{l}}L)=\mathrm{ann}_{\mathcal{C}}(\mathrm{ann}_{^{\ast }%
\mathcal{C}}(K)\ast ^{l}\mathrm{ann}_{^{\ast }\mathcal{C}}(L)).
\end{equation*}
\end{enumerate}
\end{lemma}

\begin{Beweis}
The proof of (1) is analogous to that of Lemma \ref{inn-ideal}, while (2)
follows by symmetry.$\blacksquare $
\end{Beweis}

The following result clarifies the relation between the \emph{wedge product}
and the \emph{internal coproduct} of right (left) $\mathcal{C}$-coideals
under suitable purity conditions:

\begin{proposition}
\label{wedge=}Let $A$ be a QF ring, $(\mathcal{C},\Delta ,\varepsilon )$ be
an $A$-coring and $K,L\subseteq \mathcal{C}$ be $A$-subbimodules.

\begin{enumerate}
\item Let $_{A}\mathcal{C}$ be flat and $K,L$ be right $\mathcal{C}$%
-coideals. If $_{A}L\subseteq $ $_{A}\mathcal{C}$ is pure, then $(K:_{%
\mathcal{C}^{r}}L)=K\wedge L.$

\item Let $\mathcal{C}_{A}$ be flat and $K,L$ be left $\mathcal{C}$%
-coideals. If $K_{A}\subseteq \mathcal{C}_{A}$ is pure, then $(K:_{\mathcal{C%
}^{l}}L)=K\wedge L.$

\item Let $_{A}\mathcal{C},\mathcal{C}_{A}$ be flat and $K,L\subseteq
\mathcal{C}$ be $\mathcal{C}$-bicoideals. If $_{A}K\subseteq $ $_{A}\mathcal{%
C}$ and $L_{A}\subseteq \mathcal{C}_{A}$ are pure, then%
\begin{equation}
(K:_{\mathcal{C}^{r}}L)=K\wedge L=(K:_{\mathcal{C}^{l}}L).  \label{lc=rc}
\end{equation}
\end{enumerate}
\end{proposition}

\begin{Beweis}
\begin{enumerate}
\item Assume $_{A}\mathcal{C}$ to be flat and consider the map%
\begin{equation*}
\kappa _{r}:L^{\perp (\mathcal{C}^{\ast })}\otimes _{A}K^{\perp (\mathcal{C}%
^{\ast })}\rightarrow (\mathcal{C}\otimes _{A}\mathcal{C})^{\ast },\text{ }%
[g\otimes _{A}f\mapsto (c^{\prime }\otimes _{A}c)=g(f(c^{\prime })c)].
\end{equation*}%
Then we have%
\begin{equation*}
\begin{tabular}{llll}
$(K:_{\mathcal{C}^{r}}L)$ & $=$ & $($\textrm{$ann$}$_{\mathcal{C}^{\ast
}}(K)\ast ^{r}$\textrm{$ann$}$_{\mathcal{C}^{\ast }}(L))^{\bot }$ & (Lemma %
\ref{biprod=}) \\
& $=$ & $(K^{\bot (\mathcal{C}^{\ast })}\ast ^{r}L^{\bot (\mathcal{C}^{\ast
})})^{\bot }$ & ($X,Y\subseteq \mathcal{C}$ are right coideals) \\
& $=$ & $((\Delta ^{\ast }\circ \kappa _{r})(L^{\bot (\mathcal{C}^{\ast
})}\otimes _{A}K^{\bot (\mathcal{C}^{\ast })}))^{\bot }$ &  \\
& $=$ & $\Delta ^{-1}((\kappa _{r}(L^{\bot (\mathcal{C}^{\ast })}\otimes
_{A}K^{\bot (\mathcal{C}^{\ast })}))^{\bot })$ & (\cite[Proposition 1.10
(3-c)]{Abu05}) \\
& $=$ & $\Delta ^{-1}(K^{\bot (\mathcal{C}^{\ast })\bot }\otimes _{A}%
\mathcal{C}+\mathcal{C}\otimes _{A}L^{\bot (\mathcal{C}^{\ast })\bot })$ & (
\cite[Corollary 2.9]{Abu05}) \\
& $=$ & $\Delta ^{-1}(K\otimes _{A}\mathcal{C}+\mathcal{C}\otimes _{A}L)$ & $%
A$ is cogenerator \\
& $=$ & $K\wedge L.$ &
\end{tabular}%
\end{equation*}

\item This follows from (1) by symmetry.

\item This is a combination of (1) and (2).$\blacksquare $
\end{enumerate}
\end{Beweis}

\subsection*{Fully coprime (fully cosemiprime) corings}

\qquad In addition to the notions of \emph{right} (\emph{left}) \emph{fully
coprime }and \emph{right} (\emph{left}) \emph{fully cosemiprime bicoideals},
considered as right (left) comodules in the canonical way, we present the
notion of a \emph{fully coprime }(\emph{fully cosemiprime})\emph{\ bicoideal.%
}

\begin{definition}
Let $(\mathcal{C},\Delta ,\varepsilon )$ be a non-zero $A$-coring and assume
$_{A}\mathcal{C},$ $\mathcal{C}_{A}$ to be flat. Let $0\neq B\subseteq
\mathcal{C}$ be a $\mathcal{C}$-bicomodule and consider the right $\mathcal{C%
}$-comodule $B^{r}$ and the left $\mathcal{C}$-comodule $B^{l}.$ We call $B:$

fully $\mathcal{C}$-coprime (fully $\mathcal{C}$-cosemiprime), iff both $%
B^{r}$ and $B^{l}$ are fully $\mathcal{C}$-coprime (fully $\mathcal{C}$%
-cosemiprime);

\emph{fully coprime} (\emph{fully cosemiprime}), iff both $B^{r}$ and $B^{l}$
are fully coprime (fully cosemiprime).
\end{definition}

\subsection*{The fully coprime coradical}

\qquad The \emph{prime spectra }and the associated \emph{prime radicals} for
rings play an important role in the study of structure of rings. Dually, we
define the \emph{fully coprime spectra }and the \emph{fully coprime
coradicals} for corings.

\begin{definition}
Let $(\mathcal{C},\Delta ,\varepsilon )$ be a non-zero ring and assume $_{A}%
\mathcal{C}$ to be flat. We define the \emph{fully coprime spectrum} of $%
\mathcal{C}^{r}$ as%
\begin{equation*}
\mathrm{CPSpec}(\mathcal{C}^{r}):=\{0\neq B\mid B^{r}\subseteq \mathcal{C}%
^{r}\text{ is a fully }\mathcal{C}\text{-coprime}\}
\end{equation*}%
and the \emph{fully coprime coradical }of $\mathcal{C}^{r}$ as%
\begin{equation*}
\mathrm{CPcorad}(\mathcal{C}^{r}):=\sum_{B\in \mathrm{CPSpec}(\mathcal{C}%
^{r})}B.
\end{equation*}%
Moreover, we set%
\begin{equation*}
\mathrm{CSP}(\mathcal{C}^{r}):=\{0\neq B\mid B^{r}\subseteq \mathcal{C}^{r}%
\text{ is a fully }\mathcal{C}\text{-cosemiprime}\}.
\end{equation*}%
In case $\mathcal{C}_{A}$ is flat, one defines analogously $\mathrm{CPSpec}(%
\mathcal{C}^{l}),$ $\mathrm{CPcorad}(\mathcal{C}^{l})$ and $\mathrm{CSP}(%
\mathcal{C}^{l}).$
\end{definition}

As a direct consequence of Remark \ref{coprime-tau} we get:

\begin{theorem}
\label{C-coprime-tau}Let $A$ be a QF ring and $_{A}\mathcal{C}$ be flat.
Then $\mathcal{C}^{\ast }$ is prime \emph{(}semiprime\emph{)} if and only if
$\mathcal{C}^{r}$ is fully coprime \emph{(}fully cosemiprime\emph{)}.
\end{theorem}

The following result shows that fully coprime spectrum (fully coprime
coradical) of corings is invariant under isomorphisms of corings. The proof
is analogous to that of Proposition \ref{th-coprime}.

\begin{proposition}
\label{th-sun-coprime}Let $\theta :\mathcal{C}\rightarrow \mathcal{D}$ be an
isomorphism of $A$-corings and assume $_{A}\mathcal{C},$ $_{A}\mathcal{D}$
to be flat. Then we have bijections%
\begin{equation*}
\mathrm{CPSpec}(\mathcal{C}^{r})\longleftrightarrow \mathrm{CPSpec}(\mathcal{%
D}^{r})\text{ and }\mathrm{CSP}(\mathcal{C}^{r})\longleftrightarrow \mathrm{%
CSP}(\mathcal{D}^{r}).
\end{equation*}%
In particular, $\theta (\mathrm{CPcorad}(\mathcal{C}^{r}))=\mathrm{CPcorad}(%
\mathcal{D}^{r}).$
\end{proposition}

\begin{remark}
If $\theta :\mathcal{C}\rightarrow \mathcal{D}$ is a morphism of $A$%
-corings, then it is NOT evident that $\theta $ maps fully $\mathcal{C}$%
-coprime (fully $\mathcal{C}$-cosemiprime) $\mathcal{C}$-bicoideals into
fully $\mathcal{D}$-coprime (fully $\mathcal{D}$-cosemiprime) $\mathcal{D}$%
-bicoideals, contrary to what was mentioned in \cite[Theorem 2.4(i)]{NT01}.
\end{remark}

The following example, given by Chen Hui-Xiang in his review of \cite{NT01}
(Zbl 1012.16041), shows moreover that a homomorphic image of a fully coprime
coalgebra need not be fully coprime:

\begin{ex}
Let $A:=M_{n}(F)$ be the algebra of all $n\times n$ matrices over a field $%
F, $ $B:=T_{n}(F)$ be the subalgebra of upper-triangular $n\times n$
matrices over $F$ where $n>1.$ Consider the dual coalgebras $A^{\ast
},B^{\ast }.$ The embedding of $F$-algebras $\iota :B\hookrightarrow A$
induces a surjective map of $F$-coalgebras $A^{\ast }\overset{\iota ^{\ast }}%
{\longrightarrow }B^{\ast }\longrightarrow 0.$ However, $A$ is prime while $%
B $ is not, i.e. $A^{\ast }$ is a fully coprime $F$-coalgebra, while $%
B^{\ast } $ is not (see Theorem \ref{C-coprime-tau}).
\end{ex}

\qquad As a direct consequence of Proposition \ref{Prad=CPcorad} we have

\begin{proposition}
\label{C-Prad=CPcorad}Let $A$ be a QF ring. If $_{A}\mathcal{C}$ is flat and
$\mathcal{C}^{\ast }$ is right Noetherian, then
\begin{equation*}
\mathrm{Prad}(\mathcal{C}^{\ast })=\mathrm{CPcorad}(\mathcal{C}^{r})^{\bot (%
\mathcal{C}^{\ast })}\text{ and }\mathrm{CPcorad}(\mathcal{C}^{r})=\mathrm{%
Prad}(\mathcal{C}^{\ast })^{\bot (\mathcal{C})}.
\end{equation*}
\end{proposition}

\qquad Making use of Proposition \ref{C-Prad=CPcorad}, a similar proof to
that of Corollary \ref{ME-semi} yields:

\begin{corollary}
Let $A$ be a QF ring. If $_{A}\mathcal{C}$ is flat and $\mathcal{C}^{\ast }$
is Noetherian, then%
\begin{equation*}
\mathcal{C}^{r}\text{ is fully cosemiprime}\Leftrightarrow \mathcal{C}=%
\mathrm{CPcorad}(\mathcal{C}^{r}).
\end{equation*}
\end{corollary}

\begin{corollary}
Let $A$ be a QF ring. If $_{A}\mathcal{C}$ is \emph{(}locally\emph{)}
projective and $\mathcal{C}^{r}$ is Artinian \emph{(}e.g. $\mathcal{C}_{A}$
is finitely generated\emph{)}, then

\begin{enumerate}
\item $\mathrm{Prad}(\mathcal{C}^{\ast })=\mathrm{CPcorad}(\mathcal{C}%
^{r})^{\bot (\mathcal{C}^{\ast })}$ and $\mathrm{CPcorad}(\mathcal{C}^{r})=%
\mathrm{Prad}(\mathcal{C}^{\ast })^{\bot (\mathcal{C})}.$

\item $\mathcal{C}^{r}$ is fully cosemiprime if and only if $\mathcal{C}=%
\mathrm{CPcorad}(\mathcal{C}^{r}).$
\end{enumerate}
\end{corollary}

\subsection*{Corings with Artinian dual rings}

\qquad For corings over QF ground rings several primeness and coprimeness
properties become equivalent. As a direct consequence Theorems \ref%
{E-Artinian}, \ref{C-coprime-tau} and \cite[Theorem 2.9, Corollary 2.10]{FR}
we get the following characterizations of fully coprime locally projective
corings over QF ground rings:

\begin{theorem}
\label{A-C*-Art}Let $A$ be a QF ring and $_{A}\mathcal{C},$ $\mathcal{C}_{A}$
be projective and assume $\mathcal{C}^{\ast }$ is right Artinian and $^{\ast
}\mathcal{C}$ is left Artinian. Then the following statements are equivalent:

\begin{enumerate}
\item $\mathcal{C}^{\ast }$ \emph{(}or $^{\ast }\mathcal{C}$\emph{)} is
prime;

\item $\mathcal{C}_{\mathcal{C}^{\ast }}$ \emph{(}or $_{^{\ast }\mathcal{C}}%
\mathcal{C}$\emph{)} is diprime;

\item $\mathcal{C}_{\mathcal{C}^{\ast }}$ \emph{(}or $_{^{\ast }\mathcal{C}}%
\mathcal{C}$\emph{)} is prime;

\item $\mathcal{C}^{\ast }$ \emph{(}or $^{\ast }\mathcal{C}$\emph{)} is
simple Artinian;

\item $\mathcal{C}_{\mathcal{C}^{\ast }}$ \emph{(}or $_{^{\ast }\mathcal{C}}%
\mathcal{C}$\emph{)}\ is strongly prime;

\item $\mathcal{C}^{r}$ \emph{(}or $\mathcal{C}^{l}$\emph{)} is fully
coprime;

\item $\mathcal{C}$ has non-trivial fully invariant right \emph{(}left\emph{)%
}\ $\mathcal{C}$-coideals;

\item $\mathcal{C}$ is simple.
\end{enumerate}
\end{theorem}

\qquad As a direct consequence of Theorem \ref{coprime-irr} we get

\begin{theorem}
\label{cocomm-irr}Let $C$ be a locally projective cocommutative $R$%
-coalgebra and assume $C$ to be self-injective self-cogenerator in $\mathbb{M%
}^{\mathcal{C}}.$ If $C$ is fully coprime, then $C$ is irreducible.
\end{theorem}

\subsection*{Examples and Counterexamples}

In what follows we give some examples of \emph{fully coprime }corings
(coalgebras) over arbitrary (commutative) ground rings. An important class
of fully coprime path coalgebras over fields is considered by Prof. Jara et.
al. in \cite{JMR}. For other examples of fully coprime coalgebras over
fields, the reader is referred to \cite{NT01}.

We begin with a counterexample to a conjecture in \cite{NT01}, communicated
to the author by Ch. Lomp, which shows that the converse of Theorem \ref%
{cocomm-irr} is not true in general:

\begin{c-ex}
\label{lomp}Let $C$ be a $\mathbb{C}$-vector space spanned by $g$ and an
infinite family of elements $\{x_{\lambda }\}_{\Lambda }$ where $\Lambda $
is a non-empty set. Define a coalgebra structure on $C$ by%
\begin{equation}
\begin{tabular}{lll}
$\Delta (g)=g\otimes g,$ &  & $\varepsilon (g)=1;$ \\
$\Delta (x_{\lambda })=g\otimes x_{\lambda }+x_{\lambda }\otimes g,$ &  & $%
\varepsilon (x_{\lambda })=0.$%
\end{tabular}%
\end{equation}%
Then $C$ is a cocommutative coalgebra with unique simple (1-dimensional)
subcoalgebra $C_{0}=\mathbb{C}g.$ Let $V(\Lambda )$ be the $\mathbb{C}$%
-vector space of families $\{b_{\lambda }\}_{\Lambda },$ where $b_{\lambda
}\in \mathbb{C}$ and consider the trivial extension%
\begin{equation}
\mathbb{C}\ltimes V(\Lambda )=\left\{ \left(
\begin{array}{cc}
a & w \\
0 & a%
\end{array}%
\right) \mid a\in \mathbb{C}\text{ and }w\in V(\Lambda )\right\} ,
\end{equation}%
which is a ring under the ordinary matrix multiplication and addition. Then
there exists a ring isomorphism%
\begin{equation}
C^{\ast }\simeq \mathbb{C}\ltimes V(\Lambda ),\text{ }f\mapsto \left(
\begin{array}{cc}
f(g) & (f(x_{\lambda }))_{\Lambda } \\
0 & f(g)%
\end{array}%
\right) \text{ for all }f\in C^{\ast }.
\end{equation}%
Since%
\begin{equation}
\mathrm{Jac}(C^{\ast })\simeq \mathrm{Jac}(\mathbb{C}\ltimes V(\Lambda
))=\left(
\begin{array}{cc}
0 & V(\Lambda ) \\
0 & 0%
\end{array}%
\right) ,
\end{equation}%
we have $(\mathrm{Jac}(C^{\ast }))^{2}=0,$ which means that $C^{\ast }$ is
not semiprime. So $C$ is an infinite dimensional irreducible cocommutative
coalgebra, which is not fully coprime (even not fully cosemiprime).
\end{c-ex}

\begin{punto}
(\textbf{The comatrix coring }\cite{EG-T03})\ Let $A,B$ be $R$-algebras, $Q$
a $(B,A)$-bimodule and assume $Q_{A}$ to be finitely generated projective
with dual basis $\{(e_{i},\pi _{i})\}_{i=1}^{n}\subset Q\times Q^{\ast }.$
By \cite{EG-T03}, $\mathcal{C}:=Q^{\ast }\otimes Q$ is an $A$-coring (called
the \emph{comatrix coring}) with coproduct and counit given by
\begin{equation*}
\Delta _{\mathcal{C}}(f\otimes _{B}q):=\sum_{i=1}^{n}(f\otimes
_{B}e_{i})\otimes _{A}(\pi _{i}\otimes _{B}q)\text{ and }\varepsilon _{%
\mathcal{C}}(f\otimes _{B}q):=f(q).
\end{equation*}%
Notice that we have $R$-algebra isomorphisms%
\begin{equation*}
\mathcal{C}^{\ast }:=\mathrm{Hom}_{-A}(Q^{\ast }\otimes _{B}Q,A)\simeq
\mathrm{Hom}_{-B}(Q^{\ast },\mathrm{Hom}_{-A}(Q,A))=\mathrm{End}%
_{-B}(Q^{\ast });
\end{equation*}%
and%
\begin{equation*}
^{\ast }\mathcal{C}:=\mathrm{Hom}_{A-}(Q^{\ast }\otimes _{B}Q,A)\simeq
\mathrm{Hom}_{B-}(Q,\mathrm{Hom}_{A-}(Q^{\ast },A))^{op}\simeq \mathrm{End}%
_{B-}(Q)^{op}.
\end{equation*}%
\qquad
\end{punto}

\begin{ex}
\label{Mn(A) coprime}Consider the $(A,A)$-bimodule $Q=A^{n}$ and the
corresponding comatrix $A$-coring $\mathcal{C}:=Q^{\ast }\otimes _{A}Q$
(called also the \emph{matrix coalgebra} in case $A=R,$ a commutative ring).
Then we have isomorphisms of rings%
\begin{equation*}
\mathcal{C}^{\ast }\simeq \mathrm{End}_{-A}((A^{n})^{\ast })\simeq \mathrm{%
End}_{-A}((A^{\ast })^{n})\simeq \mathbb{M}_{n}(\mathrm{End}_{-A}(A^{\ast
}))\simeq \mathbb{M}_{n}(\mathrm{End}_{-A}(A))\simeq \mathbb{M}_{n}(A),
\end{equation*}%
and%
\begin{equation*}
^{\ast }\mathcal{C}\simeq \mathrm{End}_{A-}(A^{n})^{op}\simeq \mathbb{M}_{n}(%
\mathrm{End}_{A-}(A))^{op}\simeq \mathbb{M}_{n}(A^{op})^{op}.
\end{equation*}%
Let $A$ be prime. Then $\mathcal{C}^{\ast }\simeq \mathbb{M}_{n}(A)$ and $%
^{\ast }\mathcal{C}\simeq \mathbb{M}_{n}(A^{op})^{op}$ are prime (e.g. \cite[%
Proposition 13.2]{AF74}). If moreover $A_{A}$ ($_{A}A$) is a cogenerator,
then $\mathcal{C}^{r}$ ($\mathcal{C}^{l}$) is self-cogenerator and it
follows by Remark \ref{coprime-tau} that $\mathcal{C}^{r}$ is fully coprime (%
$\mathcal{C}^{l}$ is fully coprime).
\end{ex}

\begin{ex}
Let $A\rightarrow B$ be a ring homomorphism and assume $B_{A}$ to finitely
generated and projective. Then the $A$-comatrix coring $\mathcal{C}:=B^{\ast
}\otimes _{B}B\simeq B^{\ast },$ is called the \emph{dual }$A$\emph{-coring}
of the $A$-ring $B$ as its coring structure can also be obtained from the
the $A$-ring structure of $B$ (see \cite[3.7.]{Swe75}). If $B$ is a prime
ring, then $^{\ast }\mathcal{C}:=$ $^{\ast }(B^{\ast })\simeq B$ is prime.
If moreover, $B_{A}^{\ast }$ is flat and self-cogenerator, then it follows
by analogy to Remark \ref{coprime-tau} that $^{l}\mathcal{C}$ is fully
coprime.
\end{ex}

\begin{ex}
Let $R$ be Noetherian, $A$ a non-zero $R$-algebra for which the \emph{finite
dual }$A^{\circ }\subset R^{A}$ is a pure submodule (e.g. $R$ is a Dedekind
domain) and assume $_{R}A^{\circ }$ to be a self-cogenerator. By \cite%
{AG-TW2000}, $A^{\circ }$ is an $R$-coalgebra. If the $R$-algebra $A^{\circ
\ast }$ is prime, then $A^{\circ }$ is a fully coprime $R$-coalgebra. If $A$
is a reflexive $R$-algebra (i.e. $A\simeq A^{\circ \ast }$ canonically),
then $A$ is prime if and only if $A^{\circ }$ is fully coprime.
\end{ex}

\begin{ex}
Let $A$ be a prime $R$-algebra and assume $_{R}A$ to be finitely generated
projective. Then $C:=A^{\ast }$ is an $R$-coalgebra (with no assumption on
the commutative ground ring $R$)\ and $C^{\ast }:=A^{\ast \ast }\simeq A.$
If $_{R}A^{\ast }$ is self-cogenerator (e.g. $_{R}R$ is a cogenerator), then
$A$ is a prime $R$-algebra if and only if $C$ is a fully coprime $R$%
-coalgebra.
\end{ex}

\begin{ex}
Let $R$ be a integral domain and $C:=R[x]$ be the $R$-coalgebra with
coproduct and counit defined on the generators by%
\begin{equation*}
\Delta (x^{n}):=\dsum_{j=0}^{n}x^{j}\otimes _{R}x^{n-j}\text{ and }%
\varepsilon (x^{n}):=\delta _{n,0}\text{ for all }n\geq 0.
\end{equation*}%
Then $C^{\ast }\simeq R[[x]],$ the power series ring, is an integral domain.
If moreover, $_{R}C$ is self-cogenerator (e.g. $_{R}R$ is a cogenerator),
then $C$ is fully coprime.
\end{ex}

\textbf{Acknowledgments: }Many thanks for the \emph{referee} for the careful
reading of the paper and for his/her suggestions, which improved the paper.
The author thanks also Prof. Stefan Caenepeel for inviting him to submit a
paper to this special issue of the journal. The author is grateful to Prof.
Robert Wisbauer and to Christian Lomp for fruitful discussions on the
subject. He thanks also Virginia Rodrigues for sending him preprints of \cite%
{FR}. The author is also grateful for the financial support and the
excellent research facilities provided by KFUPM.

\end{document}